\documentclass[english]{extarticle}
\usepackage[T1]{fontenc}
\usepackage{geometry}
\usepackage{color}
\usepackage{float}
\usepackage{amsmath}
\usepackage{amsthm}
\usepackage{amssymb}
\usepackage{stmaryrd}
\usepackage{esint}
\PassOptionsToPackage{normalem}{ulem}
\usepackage{ulem}

\usepackage{graphicx, caption, subcaption}
\usepackage{graphicx}
\makeatletter

\floatstyle{ruled}
\newfloat{algorithm}{tbp}{loa}
\providecommand{\algorithmname}{Algorithm}
\floatname{algorithm}{\protect\algorithmname}

\theoremstyle{plain}
\newtheorem{thm}{\protect\theoremname}
\theoremstyle{plain}

\theoremstyle{definition}

\theoremstyle{remark}

\theoremstyle{plain}

\theoremstyle{plain}
\newtheorem{lem}[thm]{\protect\lemmaname}
\ifx\proof\undefined
\newenvironment{proof}[1][\protect\proofname]{\par
\normalfont\topsep6\p@\@plus6\p@\relax
\trivlist
\itemindent\parindent
\item[\hskip\labelsep\scshape #1]\ignorespaces
}{%
\endtrivlist\@endpefalse
}
\providecommand{\proofname}{Proof}
\fi
\theoremstyle{remark}

\makeatother

\usepackage{babel}
\providecommand{\claimname}{Claim}
\providecommand{\definitionname}{Definition}
\providecommand{\lemmaname}{Lemma}
\providecommand{\propositionname}{Proposition}
\providecommand{\remarkname}{Remark}
\providecommand{\theoremname}{Theorem}
\providecommand{\corollaryname}{Corollary}
\newcommand{\x}{\textbf{x}}
\newcommand{\y}{\textbf{y}}

\newcommand{\g}{\textbf{g}}
\newcommand{\xx}{\tilde{x}}
\newcommand{\yy}{\tilde{y}}

\newcommand{\A}{\mathcal{A}}

\newcommand{\R}{\mathcal{R}}

\newcommand{\D}{\mathcal{D}}

\newcommand{\V}{\mathcal{V}}
\newcommand{\N}{\mathcal{N}}
\newcommand{\M}{\mathcal{M}}
\newcommand{\T}{\mathcal{T}}
\newcommand{\B}{\mathcal{B}}

\newcommand{\SA}{\mathcal{S}}
\newcommand{\PA}{\mathcal{P}}

\newcommand{\<}{\langle}
\newcommand{\?}{\rangle}
\newcommand{\EX}{\mathbb{E}}

\newcommand{\X}{\mathcal{X}}
\newcommand{\w}{\textbf{w}}
\newcommand{\gd}{\tilde{g}}
\newcommand{\ff}{\tilde{f}}
\newcommand{\vv}{\textbf{v}}
\newcommand{\E}{\text{Exp}}
\newcommand{\EXX}{\tilde{\mathbb{E}}_k}
\newcommand{\f}{f}
\usepackage{latexsym}
\begin{document}
\begin{center}
\textbf{\large   Distributed Optimization on Riemannian Manifolds}
\end{center}

\ \\

\begin{center}
SUHAIL MOHMAD SHAH \\
Department of Electrical Engineering,\\
Indian Institute of Technology Bombay,\\
Powai, Mumbai 400076, India.\\
(suhailshah2005@gmail.com).
\end{center}

\ \\

\noindent \textbf{Abstract} We consider the problem of consensual distributed optimization on Riemannian manifolds.  Specifically, the problem of minimization of a sum of functions subject to non-linear constraints  is studied where each function in the sum is defined on a manifold and is  only accessible at a node of a network. This generalizes the framework of distributed optimization to a much larger class of problem which subsumes the Euclidean case. Leveraging concepts from Riemannian optimization and consensus algorithms, we propose an algorithm to efficiently solve this problem. A detailed convergence analysis is carried out for (geodesically) convex optimization problems. The empirical performance of the algorithm is demonstrated using standard problems which fit the Riemannian optimization framework.\\

\noindent \textbf{Key words} Riemannian Optimization, Distributed Optimization, Manifolds, PCA.

\newpage
\section{ Introduction }

In the last few decades there has been a major effort to develop algorithms dealing with distributed optimization problems and consensus protocols. The majority of these aim at solving problems arising in wireless and sensor networks, machine learning, multi-vehicle coordination and internet transmission protocols. In such scenarios, the networks are typically spatially distributed over a large area and potentially have a large number of nodes. The absence of a central node for access to the complete system information lends them the name ``decentralized or distributed algorithms". Such a distinction is made because the lack of a central entity makes the conventional centralized optimization techniques inapplicable. This has initiated the development of distributed computational models and algorithms to support efficient operations over such networks. We refer the interested reader to \cite{nedicbook} for a succinct account of the theory of distributed optimization algorithms and the recent developments therein.

The problem of consensual distributed optimization for the Euclidean case be formally stated as :
\begin{align*}
\text{minimize} & \qquad \frac{1}{n} \sum_{i=1}^{n} f^i(x^i),\\
\text{subject to}  & \qquad \sum_{j=1}^n q_{ij}\|x^i -x^j\|^2 =0\,\,\,\forall i\in \{1,...,n\}, 
\end{align*}
where $x^i \in \mathbb{R}^n$, $f^i \,:\, \mathbb{R}^n \to \mathbb{R}$ and scalars $q_{ij}$ with $q_{ij}>0$ if the nodes $i$ and $j$ can communicate and equal to zero if not. There are usually some  additional assumptions on the network which we do not state here. As a motivating example, one can consider machine learning optimization, where the sum of functions form can arise from considering the expected risk:
\begin{equation*}
L(w) = \mathbb{E} [\ell (h(x;w),y)],
\end{equation*}
where $\ell(\cdot,\cdot) $ is the loss function so that given the input-output pair $(x, y)$, the loss incurred is $\ell(h(x;w),y) $ when $h(x;  w)$ and $y $ are the
predicted and true outputs, respectively. Lacking an access to the actual distribution on the data, we can instead consider the minimization of the empirical risk on the available data :
\begin{equation}\label{Mar-n-1}
L(w) = \frac{1}{n}\sum_{i=1}^n \ell(h(x_i;w),y_i),
\end{equation}
Setting $f^i(w) := \ell(h(x_i;w),y_i)$, we get the sum of functions objective. In solving modern machine learning problem, one has to frequently deal with distributed data sets. This happens when the data set $\{x_i,y_i\}_{i=1}^n$ may be located in disjoint subsets at different nodes of a network with only local communication. A notable example where this occurs is federated learning \cite{kone}. Thus, a learner may lack access to the entire dataset which will necessitate the deployment of distributed optimization algorithms to get a global solution.

Another  subject of interest that has received a lot of attention during the last decade is Riemannian optimization. Broadly speaking, Riemannian optimization deals with optimization problems wherein the constraint set is a Riemannian manifold. Before the advent of techniques related to this field, when dealing with such constraints, the usual practice was to alternate between optimizing in the ambient Euclidean space and projecting onto the manifold. Some classic examples include the power iteration and Oja's algorithm \cite{oja}. Both can be studied within the framework of projected gradient gradient algorithms. The potential pitfall for these schemes is that for certain manifolds (e.g., positive definite matrices), projections can be  too expensive to compute. To address these issues, there has been a lot of effort towards developing algorithms that are projection free. Riemannian optimization provides a viable alternative by directly operating on the manifold structure and exploiting the geometric structure inherent in the problem. This allows Riemannian optimization to turn the constrained optimization problem into an unconstrained one defined on the manifold. 

Having discussed the relevance of distributed optimization and Riemannian optimizaiton,  one is lead to the consideration of developing distributed optimization techniques for Riemannian manifolds. Given the plethora of practical examples which can be posed within the Rieamnnian framework such as regression on manifolds, principal component analysis (PCA), matrix completion, computing the leading eigenvector of distributed data, such an undertaking would have both practical and theoretical implications. Unfortunately, there exists a significant gap in the existing literature with regards to understanding distributed computing on Riemannian manifolds. We remark that this is a much harder problem than its Euclidean counterpart since one has to account for curvature effects. This paper aims at addressing this and begin the program of analyzing distributed optimization on Riemannain manifolds.

Formally, the problem of consensual distributed optimization for Riemannian manifolds can be posed as:

\begin{equation}\label{mainprob}
\begin{aligned}
\text{minimize} & \qquad \frac{1}{n} \sum_{i=1}^{n} f^i(x^i), \\
\text{subject to}  & \qquad x^i \in \M, \\
 & \qquad \sum_{j =1}^n q_{ij} d^2(x^i,x^j) =0\,\,\,\forall i\in \{1,...,n\},  
\end{aligned}
\end{equation}
where $\M$ is a Riemannian Manifold, $f^i \,:\, \M \to \mathbb{R}$ and $d^2(x^i,x^j)$ is the squared distance between $x^i$ and $x^j$ (see section 2 for a definition of distance for Riemannian manifolds).\\ 

\textit{Overview of the present work }: We generalize the distributed optimziation frameowrk to Riemannian manifolds and propose an algorithm for which we provide convergence guarantees. The present program of extending distributed optimization to Riemannian manifolds faces two (surmountable) hurdles :\\ 

(i) The Riemannian distance function, which is the analog of the Euclidean distance for manifolds, is not convex. This introduces local minima in the cost function that is being minimized in order to achieve consensus. These local minima constitute non-consensus configurations, by which  we mean that the limiting values at all the nodes need not be equal. \\

(ii) The gradient vectors of the cost function evaluated at different nodes are tangent vectors lying in different tangent planes. This renders the usual vector operations on them useless. This issue is addressed using the concept of parallel transport.\\

The convergence properties of the proposed algorithm are studied and convergence is established under certain assumptions.  \\

\textit{Relevant Literature }: 
\cite{Absil} traces the key ideas of gradient descent algorithms with line search on manifolds to \cite{Lue}. As computationally viable alternatives to calculating geodesics became available, the literature dealing with optimization on Riemannian manifolds grew significantly. In the last few decades, these algorithms have been increasingly challenging various mainstream algorithms for the state of art in many popular applications. As a result, there has been considerable interest in optimization and related algorithms on Riemannian manifolds, particularly on matrix manifolds. \cite{Absil} and \cite{Helmke} serve as standard references for the topic alongwith \cite{Udr} for convex optimization on Riemmanian manifolds. The iteration complexity anaysis for first order methods on Alexandrov spaces (manifolds with lower bounded curvature) was carried out in \cite{zhang}. We refer the reader to the "notes and reference" section in \cite{Absil} for a wonderful account of the history and development of this field.

The literature on distributed optimization is vast building upon the works of \cite{Tsit} for the unconstrained case and \cite{Ned} for the constrained case. \cite{Ber} was one of the earliest works to address the issue of achieving a consensus solution to an optimization problem in a network of computational agents. In \cite{Ned}, the same problem was considered subject to constraints on the  solution. The works \cite{Ram} and \cite{Kun} extended this framework and studied different variants and extensions (asynchrony, noisy links, varying communcation graphs, etc.). The non-convex version was considered in \cite{Bia}.
 
Finally, we briefly mention some of the works which study distributed consensus algorithms on Riemannian manifolds. Some of the direct applications include distributed pose estimation, camera sensor network utilization among others. \cite{Tron} lists the earliest works on Riemannian consensus algorithms as \cite{Olfati} and \cite{Scardovi}. The former considers the spherical manifold while the latter deals with the $N$-torus. \cite{Sarl} considers a more general class of compact homogeneous manifolds. However, these deal with only embedded sub-manifolds of Euclidean space. In particular they depend on the specific embedding used and hence are extrinsic. Some other relevant works which deal with coordination on Lie groups include \cite{Igar} and \cite{Sarl2}. The most important work for our purposes is \cite{Tron}. In it, the authors consider the consensus problem as an optimization problem and employ the standard Riemannian gradient descent to solve it. A primary advantage of using this approach is that it is intrinsic and hence is independent of the embedding map used.   

The paper is organized as follows : In Section 2 we provide basic definitions and recall concepts of Riemannian manifolds which we will be using here. Also, a quick review of Riemannian consensus and distributed optimization on Euclidean spaces is provided. Section 3 presents the algorithm and a brief discussion of the steps involved. Section 4 details the convergence analysis. In section 6 we test the proposed algorithm on several different examples.   \\

\textit{Some Notation }: 
Since we are dealing with distributed computation, we use a stacked
vector notation. In particular a boldface $\textbf{w}_{k}=(w_{k}^{1},\cdot\cdot\cdot ,w_{k}^{n}) \in \M^n$, where $w_{k}^{i}$ is the value stored at node $i$ during iteration $k$, with $\M^n := \M \times \cdots \times\M$ denoting the $n$-fold Cartesian product of $\M$ with itself. Additionally, for any two  tangent vectors $\textbf{a}, \textbf{b}  \in \T_\mathbf{w} \M^n$, we use the usual product metric $\langle \textbf{a} , \textbf{b}  \rangle : =  \sum_{i=1}^n \langle a^i,b^i \rangle$ under the natural identification $\T_{\mathbf{w}}\M^n =\T_{w^1}\M \bigoplus\cdots \bigoplus\T_{w^n}\M  $, which metrizes the product topology. `$\bigoplus$' denotes the Whitney sum here. This implies, geodesics, exponential maps, and gradients can be easily obtained by using the respective definitions on each copy of $\M$ in $\M^n$.

\section{Additional Background }

In this section we briefly recall the basic definitions and concepts regarding  Riemannian manifolds while also establishing notation. For a detailed treatment we refer the reader to \cite{Absil}, \cite{dcarmo} or \cite{Gall}. We also cover the basics of Riemannian consensus and Euclidean distributed optimization. 
\subsection{Basic Definitions }

\textbf{(i) Riemannian Manifolds :} Throughout this paper we let $\M$ denote a connected Riemannian manifold. A smooth $n$-dimensional manifold is a pair $(\M,\A)$, where $\M$ is a Hausdorff second countable topological space and $\A$ is a collection of charts $\{\mathcal{U}_\alpha,\psi_{\alpha}\}$ of the set $\M$ (for more details see \cite{Tron}). A Riemannian manifold is a manifold whose tangent spaces are endowed with a smoothly varying inner product $\langle \cdot,\cdot \rangle_x$ called the Riemannian metric.
The tangent space at any point $x \in \M$ is denoted by $\T_x \M$ (for definition see \cite{Absil} or \cite{dcarmo}). We recall that a tangent space admits a structure of a vector space and for a Riemannian manifold, it is a normed  vector space. The tangent bundle $\T\M$ is defined to be the disjoint union $\cup_{x \in \M}\{x\}\times\T_x\M$. The normal space at the point $x$ denoted by $\mathcal{N_{M}}(x)$ is the set of all vectors orthogonal (w.r.t to $\langle \cdot,\cdot \rangle_x$) to the tangent space at $x$. Using the norm, one can also define the arc length of a curve $\gamma : [a,b] \to \M$ as
$$L(\gamma)= \int_a^b \sqrt{\langle \dot{\gamma}(t),\dot{\gamma}(t) \rangle_{\gamma(t)}} \, dt.$$
We let $d(\cdot,\cdot)$ denote the Riemannian distance  between any two points $x,y\in \M$, i.e., 
$$ d \,:\, \mathcal{M} \times \mathcal{M} \to \mathbb{R}\, : \, d(x,y) = \inf_{\Gamma}L(\gamma),$$
where $\Gamma$ is the set of all curves in $\mathcal{M}$ joining $x$ and $y$. We recall that $d(\cdot,\cdot)$ defines a\textit{ metric }on $\M$. We also assume that the sectional curvature of the manifold is bounded above by $\kappa_{max}$ and below by $\kappa_{min}$.  \\

\textbf{(ii) Geodesics and Exponential Map :} A geodesic on $\M$ is a curve that locally minimizes the arc
length (equivalently, these are the curves that satisfy $\gamma''\in\mathcal{N_{M}}(\gamma(t))$
for all $t$). The exponential of a tangent vector $u$ at $x$, denoted
by $\text{exp}_x(u)$, is defined to be $\Gamma(1,x,u)$,
where $t\to\Gamma(t,x,u)$ is the geodesic that satisfies
\[
\Gamma(0,x,u)=x\;\text{ and }\;\frac{d}{dt}\Gamma(0,x,u)\rvert_{t=0}=u.
\]
Throughout the paper we assume a geodesically complete manifold, i.e. there is always a minimal length geodesic between any two points on the manifold. Let $\B(x,r)$ denote the open geodesic ball centered at $x$ of radius $r$, i.e.
$$\B(x,r) = \{y \in \M\,:\,d(x,y) < r \}.$$
Let $\mathcal{\tilde{I}}_x$ denote the maximal open set in $\T_x\M$ on which $\E_x$ is a diffeomorphism. On the set $\mathcal{I}_x=\E_x(\mathcal{\tilde{I}}_x)$, the exponential map 
is invertible. This inverse is called the logarithm map and we denote it by $\E_x^{-1}y$ for $y \in \mathcal{I}_x$.  \\ 

\textbf{(iii) Injectivity Radius :} The radius of the maximal geodesic ball centered at $x$ entirely contained in $\mathcal{I}_x$  is called the injectivity radius at $x$ and is denoted as $\text{inj}_x \, \M$. Also, $$ \text{inj}\, \M = \inf_x \,\{{\text{inj}_x \, \M}\}.$$

\textbf{(v) Convexity Radius  :}
The convexity radius $r_c>0$ is defined as:
$$r_c = \frac{1}{2}\Big\{ \text{inj } \M,\,\frac{\pi}{\sqrt{\kappa_{max}}} \Big\}, $$ 
where, if $\kappa_{max} \leq 0$, we set $\,1/\sqrt{\kappa_{max}}= +\infty$. This quantity plays an important role in studying the convergence properties since any open ball with radius $r\leq r_c$ is convex. A subset $\X$ of $\M$ is said to be a geodesically convex set if, given any two points in $\X$, there is an unique geodesic contained within $\X$  that joins those two points. Additionally, the function $x \mapsto d^2(y,x)$ for any fixed $y$ is (strongly) convex when restricted to $\B(y,r_c)$ (page 153, Lemma 2.9, \cite{sakai}). \\ 

\textbf{(vi) Riemannian Gradient and Hessian :} Let  $f\,:\,\M\to\mathbb{R}$ be a smooth function and $v\in \T_x\M$ be a tangent vector. We can define the directional derivative of $f$ along $v$ as $\frac{d}{dt} f(\gamma(t)) \Big|_{t=0}$ using the smooth parametrized curve $\gamma(\cdot)$ such that $\gamma(0)=x$ and $\dot{\gamma}(0) = v$. The Riemannian gradient of $f$ is then defined as the unique element of $\T_x\M$ that satisfies,
$$Df (x)[v] =\< \text{grad} \,f(x),\,v\?.  $$
The Hessian of $f$ at $x$, denoted by $\text{Hess} f(x)$, is defined as the self-adjoint linear operator which satisfies
\begin{equation*}
\frac{d}{dt^2}  f(\gamma(t) ) \Big\vert_{t=0} = \langle v, \text{Hess} f(x)v \rangle. 
\end{equation*}

\textbf{(vii) Parallel/Vector Transport :} A parallel transport $\T_x^{y}\,:\, \T_x\M\to\T_y\M$ translates the vector $\xi_x \in \T_x\M$  along the geodesic to $\T_x^{y} (\xi_x) \in \T_y\M$, while preserving norm and in some sense the direction. An important property of parallel transport is that it preserves dot products.
 
\subsection{Euclidean Distributed Optimization  }\label{edo}

Suppose we have a network of $n$ nodes/agents indexed by $1, ..., n.$ We associate with each node $i$, a function $f^i : \mathbb{R}^d \to \mathbb{R}$. Let $f:\mathbb{R}^d\to\mathbb{R}$ denote
\begin{equation}\label{H-def.}
f(\cdot) := \frac{1}{n}\sum_{i=1}^n f^i(\cdot).
\end{equation}

Let the communication network be modelled by a static undirected connected graph
$\mathcal{G=}\{\mathcal{V},\mathcal{E}\}$ where $\mathcal{V}=\{1,...,n\}$
is the node set and $\mathcal{E\subset\mathcal{V}\mathcal{\times}\mathcal{V}}$
is the set of links $(i,j)$ indicating that node $j$ can send information
to node $i$. We associate with the network a non-negative weight matrix
 $Q = [[q_{ij}]]_{i,j \in \V}$ such that
 $$q_{ij} > 0 \Longleftrightarrow (i,j)\in\mathcal{E}.$$
We let $\N_i$ denote the neighbourhood of node $i$, i.e. all the nodes such that $q_{ij}>0$. In addition, for the rest of this work, we assume the following standard assumptions on the matrix $Q$ :
\begin{enumerate}
\item[(i)] (\textit{Double Stochasticity}) $\mathbf{1}^{T}Q=\mathbf{1}^{T}$
and $Q\mathbf{1}=\mathbf{1}$.

\item[(ii)] (\textit{Irreducibility and aperiodicity})We assume that the underlying graph is irreducible, i.e., there is a directed path from any  node to any other node, and  aperiodic, i.e., the g.c.d.\ of lengths of all paths from a node to itself is one. It is known that the choice of node in this definition is immaterial. This property can be guaranteed, e.g.,  by making $q_{ii}>0$ for some $i$.

\end{enumerate}
Let $\sigma_1(Q)\geq \sigma_2(Q) \geq\cdots\geq \sigma_n(Q)\geq 0$ denote the singular values of matrix $Q$. Let $\textbf{1}_{n}$ denote the vector of all ones and $\Delta_n$ the unit simplex. Then, we have for any $x \in \Delta_n$ and $t\geq 1$, the following inequality (e.g. eq. (24), \cite{duchi})) :
$$
\big\|Q^tx - \frac{\textbf{1}_{n}}{n}\big\|_1 \leq  \sqrt{n} \big\| Q^tx - \frac{\textbf{1}_{n}}{n}  \big\|_{2} \leq  \sqrt{n}\,   \sigma^t_2(Q)
$$
 In particular, if $q^t_{ij}$ denotes the $ij$'th entry of $Q^t$, we have with $x$ equal to the \textit{i}'th co-ordinate vector $e_i$, 
 \begin{equation}\label{mix-bound}
 \big| q^t_{ij} -\frac{1}{n}\big| \leq \big\|Q^t e_i- \frac{\textbf{1}_{n}}{n}\big\|_1 \leq \sqrt{n} \, \sigma^t_2(Q)  
 \end{equation}
The objective of distributed optimization is to minimize ($\ref{H-def.}$) while simultaneously achieve consensus, i.e.
\begin{align*}
\text{minimize} & \qquad \frac{1}{n} \sum_{i=1}^{n} f^i(x^i)\\
\text{subject to} &\qquad \sum_{j \in \mathcal{N}_i}q_{ij} \|x^i -x^j\|^2 =0\,\,\,\forall i.
\end{align*}
One of the earliest works which studied the above problem and proposed a solution was \cite{Ber}. More recently, \cite{Ned} considered the constrained case. Broadly speaking, any distributed optimization algorithm in general will involve the following two steps :

\begin{enumerate}
\item[(i)] {(}\textit{Consensus Step}{)} This step involves local averaging at each node and is aimed at achieving consensus, 
\begin{equation*}
v^i_k = \sum_{j\in\N_i} q_{ij}x^j_k. 
\end{equation*}
\item[(ii)] {(}\textit{Gradient Descent Step}{)} This step is the gradient descent part aimed at minimizing $f^i$ at each node :
$$x^i_{k+1} = v^i_k - a_k \nabla f^i(v^i_k) , $$
where $a_k$ is a positive scalar decaying at a suitable rate to zero. 

\end{enumerate}
The convergence properties of the above algorithm have been extremely well studied for convex as well as the non-convex case. In Section III we modify the above scheme to give an analogous version for Riemannian Manifolds.

\subsection{Consensus on Riemannian Manifolds}
The behaviour of consensus algorithms on Riemannian manifolds is, as expected, quite different from the Euclidean case owing to curvature effects. To begin with, only local convergence can be claimed, i.e. a set of measurements located at different nodes can conform to a consensus value only if they lie within small enough convex neighbourhoods. Such a restriction does not exist for Euclidean spaces. Towards proposing a formal procedure to achieve consensus, one can follow suit from the Euclidean case and minimize the following potential function (Section III, \cite{Tron}): 
$$ \varphi(\textbf{w})= \frac{1}{2} \sum_{\{i,j\}\in\mathcal{E}} d^2(w^i,w^j),$$
where $x^i,x^j$ are the measurements in question. A Riemannian gradient descent for the above function involves the following iteration:
\begin{equation}\label{uyt}
w_{k+1}^i = \E_{w^i_k}(-\epsilon\,\text{grad}_{w_k^i}\,\varphi(\textbf{w}_k)),
\end{equation}
where $\epsilon$ is an admissible time step depending upon the upper bound $\mu_{\text{max}}$ on the Hessian of the function $\varphi(\cdot)$, and 
$
\text{grad}_{w^i}\,\varphi(\textbf{w}_k) =- \sum_{j \in \N_i} \E^{-1}_{w^i} w^j.
$ The failure of achieving consensus with the above iteration occurs due to $\varphi$ having multiple critical points globally, some of which may be local minima corresponding to a non-consensual configuration. This can happen even for the simple case of a circle, where convergence to a consensus configration can take place only if the measurements lie within a semi circle \cite{Sarl12}. 

By definition, for undirected connected graphs, the global minima of $\varphi(\cdot)$ can be seen to belong to the consensus sub-manifold $\D$:
\begin{equation*}
\D:=\{(w^i,....,w^n) \in \M^n\,: \,w^i =w^j \,\forall\,i,j\}. 
\end{equation*}
This set is the diagonal space of $\M^N$ and it represents all possible consensus configurations of the network. Let us define the following the set:
\begin{equation*}
\SA = \Big\{ \textbf{x}\in \M^n\,:\ \exists y \in \M \text{ s.t. } d(y,x^i) < \frac{r_c}{2} \forall \, i \Big\},
\end{equation*} 
where $r_c$ is the convexity radius. We note that this implies:
\begin{align*}
d^2(x^i,x^j) &\leq  2d^2(x^i,y) + 2d^2(y,x^j)\\
&\leq r_c^2
\end{align*}  
so that $x^i\, \forall\, i$ lie within  a convex neighbourhood. For set $\SA$, $\varphi(\cdot)$ has only global minima  (Theorem 5, \cite{Tron}). Thus, one can reasonably expect gradient descent to converge to it, provided the iterates $x^i_k$, generated by (\ref{uyt}), stay in $\mathcal{S}$. Assuming the latter fact, one can confirm the convergence by performing  a second order Taylor expansion of the smooth map $\varphi\circ\E(\cdot)$, using the expression $w^i_{k+1} =  \E_{w^i_k}(-\epsilon\,\text{grad}_{w^i_k}\,\varphi(\textbf{w}_k))$:\footnote{The Taylor Expansion is for $\varphi(\cdot): \M^n \to \mathbb{R}$. So, $\text{grad }\,\varphi(\w)$ belongs to $\T_\w \M^n$.}
\begin{align*}
\varphi(\textbf{w}_{k+1}) &\leq \varphi(\textbf{w}_k)+ \langle \text{grad}\, \varphi(\mathbf{w}_k)),  \E_{\w_k}^{-1} (\w_{k+1})\rangle  + \frac{\mu_{\text{max}}\|\text{grad}\, \varphi(\mathbf{w}_k)) \|^2}{2}\epsilon^2\\
\varphi(\textbf{w}_{k+1}) &\leq \varphi(\textbf{w}_k)-\|\text{grad}\, \varphi(\mathbf{w}_k)) \|^2\,\epsilon   + \frac{\mu_{\text{max}}\|\text{grad}\, \varphi(\mathbf{w}_k)) \|^2}{2}\epsilon^2\\
&\leq \varphi(\textbf{w}_k) - \|\text{grad}\, \varphi(\mathbf{w}_k)) \|^2\,\epsilon\, \Big( 1 - \frac{\mu_{\text{max}} \epsilon}{2}\Big),
\end{align*}
where $\mu_\text{max}$ is an upper bound on the Hessian of $\varphi$ (see Prop. 9, \cite{Tron} for a bound on this quantity). Thus, $\varphi(\w_{k+1}) < \varphi{(\w_k)}$ as long as $\epsilon \in (0,2\mu^{-1}_{\text{max}})$. \\

\section{ A Distributed Optimization algorithm on Riemannian Manifolds}

The pseudo code of the proposed algorithm to solve (\ref{mainprob}) is detailed in Algorithm 1.  We discuss the steps in the following:

\begin{enumerate}

\item [(S1)]  This constitutes the consensus part of the algorithm. At each node $i$, the gradient of the function $\varphi(\cdot)$ is evaluated with respect to $w^i_k$. The estimate at node $i$ moves from the current estimate $w^i_k$ to a new estimate, designated as $v^i_k$, along the exponential in the direction of the negative gradient with an admissible step size $\epsilon$.

\item [(S2)]  The quantity $g_k^i$ helps to keep track of the average (Riemannian) gradient $\nabla f(\cdot)$ via a `consensus in the tangent space'. For each node $i$, this step essentially involves gathering the parallel transported estimates $\{g_k^j\}_{j \in \N(i)}$ of the neighbours and performing a weighted average of them according to the network weights $q_{ij}$. To this weighted average, we add the difference of the current gradient $\nabla f^i (v_k^i)$and the previous gradient $\nabla f^i(v_{k-1}^i)$. This technique is sometimes referred to as the gradient tracking technique (see e.g. \cite{qu}, \cite{di}, \cite{Ned2017}). We note that this step is fully distributed and only involves local computations.

\item [(S3)]  This step constitutes the optimization part of the algorithm. It is the standard Riemannian gradient descent algorithm performed using the exponential map with a suitably fast decaying time step $a_k$. The decaying time step is necessary here because of the need to achieve consensus and to provide exact convergence guarantees. For a constant time step, one can only establish convergence within a neighbourhood of the optimum.

\end{enumerate}

\begin{algorithm}[H]
\textbf{Input :}

Manifold $\M$; Cost functions $f^i(\cdot)$ for all $i$, Exponential map $\text{Exp}(\cdot)$ or Retraction $\mathcal{R}(\cdot)$; Parallel transport $\mathcal{T}(\cdot)$; Time steps $\epsilon$ and $a_{k}$; Network graph $\mathcal{G}=(\mathcal{V},\mathcal{E})$
with weight matrix $Q$.

\textbf{Initial Conditions:}

Initialize $\mathbf{w}_{0}$ with $w_0^i =w_0^j$ and $g_0^i := \nabla f^i(v_0^i) \,\forall i$.

\textbf{For $k=0,1,2.....$ do:}

At each node $i\in\{1,..,n\}$ do:\\
(S1) \textbf{{[}Riemannian Consensus Step{]} :}
$$
v_{k}^{i}=\E_{w_{k}^{i}}(-\epsilon \, \text{grad}_{w_{k}^{i}}\varphi(\mathbf{w_{k}})). 
$$
(S2) \textbf{[Gradient Update Step] :}
\begin{equation*}
g_{k}^i =\sum_{j=1}^{n}q_{ij} \T_{v_{k-1}^j}^{v_{k}^i}  \,g_{k-1}^j + \nabla f^i(v_{k}^{i}) - \T_{v_{k-1}^i}^{v_{k}^i}  \nabla f^i(v_{k-1}^{i}).
\end{equation*}
(S3) \textbf{{[}Gradient Descent Step{]} :}
$$w_{k+1}^{i}=\E_{v_{k}^{i}}(-a_{k} g^i_k).$$
$k\leftarrow k+1$

\textbf{end}

\caption{Distributed Optimization on Riemannian Manifolds }
\end{algorithm}

Alternatively, one could use a (second order) retraction map  instead of $\E(\cdot)$ in Step 1 and 3 of Algorithm 1 with essentially no change in the behaviour of the algorithm. The main advantage of a retraction is a lower computational overhead. A retraction on $\mathcal{M}$ is a smooth mapping $\mathcal{R}:\mathcal{TM}\to\mathcal{M}$, where $\mathcal{TM}$ is the tangent bundle,
with the following properties :

i) $\mathcal{R}_{x}(0_{x})=x$, where $\mathcal{R}_{x}$ is the restriction of the retraction
to $\mathcal{T}_{x}\mathcal{M}$ and $0_{x}$ denotes the zero element
of $\mathcal{T}_{x}\mathcal{M}$.

ii) With the canonical identification $\mathcal{T}_{0_{x}}\mathcal{T}_{x}\mathcal{M}\simeq\mathcal{T}_{x}\mathcal{M}$,
$\mathcal{R}_{x}$ satisfies
\[
D\mathcal{R}_{x}(0_{x})=\text{id}_{\mathcal{T}_{x}\mathcal{M}},
\]
where $\text{id}_{\mathcal{T}_{x}\mathcal{M}}$ denotes the identity
mapping on $\mathcal{T}_{x}\mathcal{M}$. The fact that a retraction mapping does not affect the behaviour of the algorithm hinges on two properties. The first of these is:
\begin{equation*}
 d(x,R_x(\xi_x)) \leq c \| \xi_x\|.
\end{equation*}
This obviously holds for the exponential with an equality for $c=1$. The second one is 
$$
d(R_x(t \xi_x)\E_x(t \xi_x) ) =\mathcal{O}(t^2),
$$
i.e., the retraction is a second order one. On a related notion, generalizing the notion of parallel transport, a vector transport associated with a retraction $\R$  is defined as a smooth mapping  ($\oplus $ denotes the Whitney sum, we refer the reader to page 169, \cite{Absil} for exact definition)
$$\T\M \oplus \T\M \to \T\M\, \,:\,( \eta_x,\xi_x) \mapsto \,\T_{\eta_x}(\xi_x) \in \T\M $$
satisfying :

i) If $\pi(\cdot)$ denotes the foot of a tangent vector, then $\pi(\R_x (\eta_x)) =  \pi(\T_{\eta_x}(\xi_x)). $,

ii) $\T_{\eta_x}\,:\,\T_{x}\M \to \T_{\R(\eta_x)}\M  $ is a linear map,

iii) $\T_{0_x}(\xi_x) = \xi_x$ 

Vector transports are in general more computationally appealing than parallel transport. The vector transport is well defined as long as we stay within the injectivity radius defined with respect to the retraction.

\section{ Convergence Analysis }
In this section, we analyze the convergence properties of Algorithm 1 for convex optimization on Rimeannian manifolds. This class of problems has found applications in wide variety of fileds (see e.g., Chapter 4, \cite{Udr}, Chapter 7, \cite{boumal}). To be more precise, this means we consider problems where each $f^i(\cdot)$ (and hence $f(\cdot)$, see Theorem 3.3, \cite{Udr}) is geodesically convex and $x \in \X\subset \M$, with $\X$ being a compact convex set. A function $f : \M \to R$ is said to be geodesically convex if for any $x, y \in \M$, a geodesic $\gamma$ such that $\gamma(0) = x$ and $\gamma(1)= y$, and $t \in [0, 1]$, the following holds:
$$
f(\gamma(t)) \leq (1 - t)f(x) + tf(y).
$$
An equivalent definition for geodesic convexity is the existence of a vector $\partial f(x) \in \T_x \M$, such that
\begin{equation}\label{conv}
f(y) \geq f(x ) + \langle \partial f(x) , \E_{x}^{-1}y \rangle. 
\end{equation}
We note that for a differentiable $f:\M \to \mathbb{R}$, we have  $\partial f (x)= \nabla f(x)$, where $\nabla f$ is the Riemannian gradient. We remark that if we are using a retraction, then the concept of  a convexity has to be defined accordingly (see Definition 3.1, \cite{Huang}).

Throughout this section we will assume that $f$ is geodesically $L_f$-smooth by which we mean, for any $x,\,y \in \M$
\begin{equation}\label{smooth}
f(y) \leq f(x) + \langle \nabla f(x), \E^{-1}_{x} (y)  \rangle + \frac{L_f}{2} \|\E^{-1}_{x} (y) \|^2
\end{equation}
This can be shown to imply:
\begin{equation}\label{papallel}
\| \nabla f(x) - \T_{y}^x \nabla f(y)\| \leq L_f d(x,y).
\end{equation}
Since $\X$ is compact, we can assume that the exponential map is invertible (Hopf-Rinow theorem). We remark that these conditions are standard assumptions for Riemannian optimization and most practical problems of interest satisfy them.

Using the stacked vector notation, one can write step (S1) and (S3) of Algorithm 1 as:
\begin{align*}
\vv_{k}&=\E_{\w_{k}}(-\epsilon \, \text{grad}\,\varphi(\mathbf{\w_{k}})) \\
\w_{k+1} &=\E_{\vv_{k}}(-a_{k} \g_k),
\end{align*}
where, 
\begin{equation*}
\text{grad}\,\varphi(\w_k) = ( \text{grad}_{w^1_k}\,\varphi(w^i_k),\cdots,\text{grad}_{w_k^n}\,\varphi(w^n_k)) \in \T_{\w_k}\M^n 
\end{equation*}
and
\begin{equation*}
 \E_{\w_k}\big(-\epsilon\text{grad} \,\varphi(\mathbf{\w_{k}}))\big) :=\\
\Big( \E_{w^1_k}\big(-\epsilon\text{grad}_{w^1_{k}}\varphi(\mathbf{\w_{k}})),\cdots, \E_{w^n_k}\big(-\epsilon\text{grad}_{w^n_{k}}\varphi(\mathbf{\w_{k}}))\big)\Big).
\end{equation*}
Since we assume $\X$ is compact, we have from the smoothness of $f$:  
\begin{equation}\label{funcbound}
\| \nabla f(\vv_k)\| \ \leq C.
\end{equation}
for some positive constant $C< \infty$. We also take $\|\g^i_k\| \leq C$. This again is a reasonable assumption given the boundedness of $\nabla f(\vv_k)$, see Lemma \ref{lem2} for a bound on the difference between $\g_k^i$ and $\nabla f(\cdot)$. 
We first prove that Algorithm 1 achieves consensus. 
\subsection{Consensus}
The following lemma shows that Algorithm 1 asymptotically achieves consensus if the time steps are properly chosen. 

\begin{lem} \label{con}
Suppose $\epsilon \in (0,2\mu^{-1}_{\text{max}})$, $\sum_{k=1}^\infty a_k =\infty$ and $\sum_{k\geq 0} a_k^2 < \infty$. Then, we have,
\begin{equation}\label{gradconv}
\lim_k \|\text{grad} \,\varphi(\mathbf{w_k})\|^2 = 0,\,\,\,\forall \,i.
\end{equation}
so that, $$ \textbf{w}_k \to \mathcal{D}\,\text{ as }\,k\to \infty,$$
 i.e.,
$$d(w_k^i,w_k^j) \to 0 \,\text{ as }\,k\to \infty,\,\forall \,i,j $$
\end{lem}

\begin{proof}
We begin by using the smoothness property for $\varphi(\cdot)$ within $\X$:
\begin{multline}\label{lem1-1} 
 \varphi(\textbf{v}_{k}) \leq \varphi(\textbf{w}_k)- \|\text{grad} \,  \varphi(\mathbf{w}_k) \|^2\,\epsilon  + \frac{\mu_{\text{max}}\|\text{grad} \, \varphi(\mathbf{w}_k) \|^2}{2}\epsilon^2
 \leq \varphi(\textbf{w}_k)-\|\text{grad} \, \varphi(\mathbf{w}_k) \|^2\,\Big( 1 -\frac{ \mu_{\text{max}} \epsilon}{2}\Big)\epsilon,
\end{multline}
where we used the fact that $\E^{-1}_{\w_k}(\vv_k) = - \epsilon\text{grad} \, \varphi(\mathbf{w}_k)$. Similarly, using $\w_{k+1} =\E_{\vv_{k}}(-a_{k} \g_k)$ gives:
\begin{align} \label{2-10}
\varphi(\w_{k+1}) &\leq \varphi(\vv_{k}) - a_k \< \text{grad} \, \varphi(\mathbf{v}_k) , \g_k \?  + \frac{\mu_{\text{max}}\| \g_k \|^2}{2} a_k^2. \nonumber \\
&\leq \varphi(\vv_{k}) +  \frac{\eta \|\text{grad} \, \varphi(\mathbf{v}_k)\|^2}{2}  + \frac{\| \g_k \|^2}{2\eta} a_k^2     + \frac{\mu_{\text{max}}\| \g_k \|^2}{2} a_k^2.
\end{align}  
where we have used the fact that $2\langle a,b \rangle \leq \eta \|a\|^2 + \frac{\|b\|^2}{\eta}  $ for any $\eta>0$. We have
\begin{equation}\label{vw}
 d^2(\w_{k},\vv_k) =  d^2 \big( \w_{k}, \E_{\w_k}(-\epsilon \text{grad} \,\varphi (\textbf{w}_{k})\big) 
 \leq   \epsilon^2 \| \text{grad}\,\varphi (\textbf{w}_{k})\|^2
\end{equation} 
 Let $\T_{\w}^{\vv}\text{grad} \, \varphi (\textbf{w}) :=  [\T_{w^1}^{v^1} \text{grad}_{w^1}\varphi (\textbf{w}),...,\T_{w^n}^{v^n}\text{grad}_{w^n} \varphi (\textbf{w})]$. The smoothness of $\varphi(\cdot)$ implies
\begin{equation}\label{v-field}
  \big\| \text{grad}\, \varphi (\textbf{v}_k) - \T_{\w_k}^{\vv_k}\text{grad} \, \varphi (\textbf{w}_k) \big\| \leq L_{\varphi}d(\vv_k,\w_k),
\end{equation}
where $L_{\varphi}$ is the Lipschitz constant of $\text{grad} \,\varphi(\cdot)$. We remark that the above inequality uses the fact that the Hessian of $\varphi$ has a bounded operator norm $\mu_\text{max}$, which implies Lipschitz continuity of $\text{grad}\, \varphi$ with Lipschitz constant $L_{\varphi} \geq \mu_{\text{max}}$ (Corollary 10.41, \cite{boumal}). Combining (\ref{vw}) and (\ref{v-field}), we have the bound,
\begin{equation}\label{newb}
 \big\| \text{grad} \, \varphi (\textbf{v}_k) - \T_{\w_k}^{\vv}\text{grad} \,\varphi (\textbf{w}_{k}) \big\|^2 \leq  \epsilon^2 L^2_{\varphi} \, \| \text{grad} \,\varphi (\textbf{w}_{k})\|^2.
\end{equation}
We use (\ref{newb}) in (\ref{2-10}) as follows:
\begin{align*}
\varphi(\w_{k+1})  &\leq \varphi(\vv_{k}) +  \frac{\eta \|\text{grad}\, \varphi(\mathbf{v}_k) -\text{grad} \, \varphi(\mathbf{w}_k)  +\text{grad} \, \varphi(\mathbf{w}_k)  \|^2}{2}   + \frac{ \| \g_k\|^2}{2\eta} a_k^2     +
 \frac{\mu_{\text{max}}\| \g_k \|^2}{2} a_k^2.\\
&\leq \varphi(\vv_{k}) +  \eta \|\text{grad} \, \varphi(\mathbf{v}_k) - \text{grad}\, \varphi(\mathbf{w}_k)  \|^2   + \eta \|  \text{grad} \,\varphi(\mathbf{w}_k)\|^2+ \Big(\frac{1}{2\eta}     + \frac{\mu_{\text{max}}}{2}\Big)C^2 a_k^2\\
&\leq  \varphi(\vv_{k}) +  \eta \big( 1 + \epsilon^2 L_{\varphi}^2  \big) \|  \text{grad}\, \varphi(\mathbf{w}_k)\|^2+ \Big(\frac{1}{2\eta}     + \frac{\mu_{\text{max}}}{2}\Big)C^2 a_k^2\\
&\leq \varphi(\w_{k}) -  \|\text{grad} \,\varphi(\mathbf{w}_k) \|^2\,\epsilon\, \Big( 1 - \frac{\mu_{\text{max}} \epsilon}{2}\Big) +  \eta \big( 1 + \epsilon^2 L_{\varphi}^2  \big) \|  \text{grad}\, \varphi(\mathbf{w}_k)\|^2 +  \Big(\frac{1}{2\eta}     + \frac{\mu_{\text{max}}}{2}\Big)C^2 a_k^2,
\end{align*}
where we have used (\ref{lem1-1}) in the last inequality. Setting $\eta = \frac{ \epsilon}{2(1 + \epsilon^2 L_{\varphi}^2)}$, we get
\begin{align*}
\varphi(\w_{k+1})   &\leq \varphi(\w_{k}) -  \frac{\epsilon}{2}\, \Big( 1 - \frac{\mu_{\text{max}}\epsilon}{2} \Big)\|  \text{grad} \, \varphi(\mathbf{w}_k)\|^2  + \Bigg(\frac{(1 + \epsilon^2 L_{\varphi}^2)}{\epsilon}     + \frac{\mu_{\text{max}}}{2}\Bigg)C^2 a_k^2.
\end{align*}
Re-arranging gives
\begin{align}\label{transient-1}
\frac{\epsilon}{2}\, \Big( 1 - \frac{\mu_{\text{max}}\epsilon}{2} \Big)\,\|  \text{grad} \, \varphi(\mathbf{w}_k)\|^2   &\leq \big( \varphi(\w_{k}) - \varphi(\w_{k+1})\big)
   + \Bigg(\frac{(1 + \epsilon^2 L_{\varphi}^2)}{\epsilon}     + \frac{\mu_{\text{max}}}{2}\Bigg)C^2 a_k^2
\end{align}
Summing the above equation from $k=1$ to $\ell$ gives:
\begin{multline*}
 \Big( 1 - \frac{\mu_{\text{max}}\epsilon}{2} \Big)\frac{\epsilon}{2} \sum_{k=0}^\ell  \|  \text{grad}\, \varphi(\mathbf{w}_k)\|^2   \leq  \big( \varphi(\w_{0}) - \varphi(\w_{\ell})\big)
+ \Bigg(\frac{(1 + \epsilon^2 L_{\varphi}^2)}{\epsilon}     + \frac{\mu_{\text{max}}}{2}\Bigg)C^2\sum_{k=0}^\ell  a_k^2\\
\leq \Bigg(\frac{(1 + \epsilon^2 L_{\varphi}^2)}{\epsilon}     + \frac{\mu_{\text{max}}}{2}\Bigg)C^2 \sum_{k=0}^\ell a_k^2,
\end{multline*}
where we used the fact that $\varphi(\w_0) = 0$ and $\varphi(\w_{l})\geq 0$ to obtain last inequality. Taking $\ell \to \infty$ gives:
\begin{align*}
 \frac{\epsilon}{2}\, &\Big( 1 - \frac{\mu_{\text{max}}\epsilon}{2} \Big) \sum_{k=0}^\infty\|  \text{grad}\, \varphi(\mathbf{w}_k)\|^2   \leq  \Bigg(\frac{(1 + \epsilon^2 L_{\varphi}^2)}{\epsilon}     + \frac{\mu_{\text{max}}}{2}\Bigg)C^2 \sum_{k=0}^\infty a_k^2 <\infty
\end{align*}
which implies $\lim_k \|  \text{grad} \, \varphi(\mathbf{w}_k)\|^2 = 0$, since $\epsilon < \frac{2}{\mu_\text{max}}$. This proves the result since $\lim_k  \|  \text{grad} \, \varphi(\mathbf{w}_k)\| = 0$ implies any limit point of $\w_k$ reaches a (global) minimum of $\varphi$. 
\end{proof}
\subsection{Geodesically Convex Function Convergence Analysis}
To prove the main convergence result, we first establish some preliminary results.

\begin{lem}\label{lem1}
The sequence $\{ \text{grad} \,\varphi(\w_k)\}_{k \geq 0}$ satisfies the following relation:
\begin{align}\label{transient-101}
\epsilon\,\|  \text{grad}\,\varphi(\mathbf{w}_k )\|^2   &\leq \frac{2}{ \Big( 1 - \frac{\mu_{\text{max} \epsilon}}{2} \Big)}\Bigg\{ \big( \varphi(\w_{k}) - \varphi(\w_{k+1})\big)
   +  \frac{2}{ \Big( 1 - \frac{\mu_{\text{max} }\epsilon}{2} \Big)} \Bigg(\frac{(1 + \epsilon^2 L_{\varphi}^2)}{\epsilon }     + \frac{\mu_{\text{max}}}{2}\Bigg)C^2a_k^2 \Bigg\}
\end{align}

\end{lem}

\begin{proof}
This is just equation (\ref{transient-1}).
\end{proof}

We recall the notation $\nabla f(x) := \frac{1}{n} \sum_{i=1}^n \nabla f^i (x)$. In the next lemma, we examine the sequence $\{g^i_k\}_{k\geq 0}$. To prove it, we assume the follwoing mild condition on the time step:
$$
\frac{a_{k}}{a_{k+1}} \leq a
$$
for some constant $a$. As an example, for $a_k = \frac{1}{k}$, $a$ can be taken equal to 2.
\begin{lem}\label{lem2}
For all $i$, the sequence $\{g^i_k\}_{k\geq 0}$ satisfies the following relation:
$$
\sum_{m=0}^k a_m \|g_m - \nabla f(v_m^i)\| \leq \frac{\sqrt{n}M}{1-\sigma_2(Q)}\sum_{m=0}^k a_m^2,
$$ 
where,
\begin{equation}\label{M}
M:= 1 +aC + \frac{2\epsilon\Big(1+\frac{9}{2\epsilon^2} \Big)}{1-\frac{\mu_{\text{max}}\epsilon}{2} } \Big( \frac{(1 + \epsilon^2 L_{\varphi}^2)}{ \epsilon}     + \frac{\mu_{\text{max}}}{2}\Big)C^2 
\end{equation}

\end{lem}

\begin{proof}
The proof is provided in Appendix \ref{A2}.

\end{proof}
The next result establishes the local convexity of $\varphi(\cdot)$. 
\begin{lem}
Let $\x, \y \in \M$ be any two points in $\M^n$ with $d(x^i,y^j) < r^*$ for all components $x^i,y^i,\,i=1$ to $n$. Then, we have 
\begin{equation*}
\varphi(\y) \geq \varphi(\x) + \langle \text{grad}_{\x} \varphi(x), \E^{-1}_{\x}(\y)\rangle,
\end{equation*}
for any $r^* \in (0,\infty)$ if the upper bound on the sectional curvature satisfies $\kappa_{\text{max}} \leq 0$ and $r^*\in \Big(0,\frac{\pi}{2\sqrt{\kappa_{\text{max}}}} \Big)$, if $\kappa_{\text{max}}> 0$.
\end{lem}
\begin{proof}
The proof is provided in Appendix II.
\end{proof}
Before stating the main result, we recall the following lemma (see e.g. Lemma 5, \cite{zhang} or  \cite{bonnabel}) which states the law of cosines for non-Eucliedean geometry:
\begin{lem}
If $a, b, c$ are the side lengths of a geodesic triangle in a Riemannian manifold and $A$ is the angle between sides $b$ and $c$ defined through inverse exponential map and inner product in tangent space, then
\begin{equation}\label{triangle}
a^2 \leq \xi b^2  -2bc \cos(A)  + c^2,
\end{equation}
where
\begin{equation}\label{xi}
\xi :=
 \frac{ \sqrt{|\kappa_{\text{min}}|} D}{\tanh (  \sqrt{|\kappa_{\text{min}}|} D)}
\end{equation}
with $\max_{x,y \in \M} d(x,y) \leq D $ and $\kappa_{\text{min}}$ being a lower bound on the sectional curvature. \\
\end{lem}
We are now in a position to prove the main result which we can use to prove the convergence of Algorithm 1. We use the notation $f(v) = n^{-1}\sum_{i=1}^n f^i(v)$ for $v\in \M$ and $f(\vv) := n^{-1}\sum_{i=1}^n f^i(v^i)$ for $v \in \M^n$.

\begin{thm}
Suppose $\epsilon \in (0,2\mu^{-1}_{\text{max}})$ and $\sum_{k\geq 0} a_k = \infty$, $\sum_{k\geq 0} a_k^2 <\infty$. Let $w^*$ denote a optimal point of (\ref{mainprob}) belonging to $\X$, where $\X \subset B(w^*,r^*)$, with $r^*$ as in Lemma 4. Then, the objective value sequence $\{f(\vv_k)\}_{m=1}^k$ for all $i$ satisfies the bound:
\begin{multline}\label{f-bound}
\sum_{m=0}^k 2 a_k \big( \f(\vv_k)- f(w^*)\big) \leq   d^2(\w_0 , \w^*)   
+ \Bigg\{  \frac{2\sqrt{n}MD}{1-\sigma_2(Q)} +  \xi C^2   \\+  \frac{2}{ \Big( 1 - \frac{\mu_{\text{max} }\epsilon}{2} \Big)}  \Bigg(\frac{(1 + \epsilon^2 L_{\varphi}^2)}{\epsilon }     + \frac{\mu_{\text{max}}}{2}\Bigg) \epsilon \xi  C^2 \Bigg\} \sum_{m=0}^k a_k^2,
\end{multline}
\end{thm}
\begin{proof}
Consider the geodesic triangle determined by $w^i_{k+1},\,v_k^i$ and $w^*$. Then, it can be seen that from $w^i_{k+1} = \E_{v_k^i}(-a_k g_k^i)$, the angle $\angle w^i_{k+1}v_k^iw^* $ satisfies:
$$  
\cos (\angle w^i_{k+1}v_k^iw^*) =  \langle-a_k \,g^i_k, \E^{-1}_{v_k^i} ( w^*) \rangle 
$$
Then, using Lemma 5 with $a= d(w^i_{k+1},w^* ) ,\, b =  \| \E^{-1}_{ v^i_{k}} (w_{k+1}^i) \|$, $c = d(v^i_k ,w^*)$, and $A= \angle w^i_{k+1}v_k^iw^* $, we have 
\begin{align*}
d^2(w^i_{k+1},w^* ) &\leq d^2(v^i_k ,w^*) +2 a_k\langle g^i_k, \E^{-1}_{v_k^i} ( w^*) \rangle + \xi \| \E^{-1}_{ v^i_{k}} (w_{k+1}^i) \|^2 \\
&\leq d^2(v^i_k ,w^*) +2a_k \langle \nabla f(v_k^i), \E_{v_k^i} ( w^*) \rangle + g^i_k - \nabla f(v_k^i), \E^{-1}_{v_k^i} ( w^*) \rangle + \xi a_k^2 \| g_k^i  \|^2 \\ 
&\leq d^2(v^i_k ,w^*) + 2a_k (f(w^*)-f(v_k^i)) + 2a_k \langle g^i_k - \nabla f(v_k^i), \E^{-1}_{v_k^i} ( w^*) \rangle ++ \xi a_k^2 \| g_k^i  \|^2 .
\end{align*}
which gives
\begin{multline*}
2a_k (f(v_k^i) -f(w^*)) \leq d^2(v^i_k ,w^*) -d^2(w^i_{k+1},w^* ) + 2a_k D \| g^i_k - \nabla f(v_k^i)\| +  \xi C^2 a_k^2
\end{multline*}
Summing over $i$ gives:
\begin{multline}\label{proof111}
2a_k (\f(\vv_k) -f(w^*)) \leq d^2(\vv_k ,\w^*) -d^2(\w_{k+1},\w^* ) + 2a_kn^{-1} D\sum_{i=1}^n \| g^i_k - \nabla f(v_k^i)\| +  \xi C^2 a_k^2,
\end{multline}
where $\f(\vv_k) : = n^{-1}\sum_{i=1}^n f^i(v_k^i)$. Similarly, with $\vv_k = \E_{\w_k} ( -  \text{grad}_{\w_{k}} \varphi(\mathbf{w}_k))$, we have:
\begin{align*}
d^2(\vv_{k},\w^* ) &\leq d^2(\w_k ,\w^*) + \epsilon  \langle  \text{grad}\, \varphi(\mathbf{w}_k), \E^{-1}_{\w_k} ( \w^*) \rangle + \xi \| \E^{-1}_{\w_k} ( \vv_{k})  \|^2\\ &\leq d^2(\w_k ,\w^*) + \epsilon \big(  \varphi(\w^*) -\varphi(\w_{k}) \big)+ \xi \| \E^{-1}_{\w_k} ( \vv_{k})  \|^2\\
&\leq d^2(\w_k ,w^*) +  \xi  \epsilon^2 \|  \text{grad}\,\varphi(\mathbf{w}_k) \|^2,
\end{align*}
where we have used Lemma 4 in the second inequality and the fact that $\varphi(\textbf{w}^*) = 0 $ with $\w^*:= (w^*,\cdots,w^*)$ and $\varphi(\textbf{w}_k)\geq 0$  for obtaining the third inequality. Using Lemma \ref{lem1}, we have
\begin{align*}
d^2(\vv_{k}, \w^* ) &\leq d^2(\w^i_k , \w^*)  + \frac{ 2 \xi \epsilon }{1 - \frac{\mu_{\text{max}} \epsilon}{2}}\Bigg( \varphi(\w_k ) - \varphi(\w_{k+1}) \Bigg)+  \frac{2}{ \Big( 1 - \frac{\mu_{\text{max} }\epsilon}{2} \Big)} \Bigg(\frac{(1 + \epsilon^2 L_{\varphi}^2)}{\epsilon }     + \frac{\mu_{\text{max}}}{2}\Bigg) 2\xi  \epsilon C^2 a_k^2
\end{align*}
Plugging this expression for $d^2(\vv_k,\w^* )$ in (\ref{proof111}) gives:
\begin{multline*}
2a_k (\f(\vv_k) -f(w^*)) \leq d^2(\vv_k ,\w^*) -d^2(\w_{k+1},\w^* ) +  2a_k n^{-1}D \sum_{i=1}^n \| g^i_k - \nabla f(v_k^i)\| +  \xi C^2 a_k^2 \\ + \frac{ 2 \xi \epsilon }{1 - \frac{\mu_{\text{max}} \epsilon}{2}}\Bigg( \varphi(\w_k ) - \varphi(\w_{k+1}) \Bigg)    +  \frac{2}{ \Big( 1 - \frac{\mu_{\text{max} }\epsilon}{2} \Big)} \Bigg(\frac{(1 + \epsilon^2 L_{\varphi}^2)}{\epsilon }     +  \frac{2}{ \Big( 1 - \frac{\mu_{\text{max} }\epsilon}{2} \Big)} \frac{\mu_{\text{max}}}{2}\Bigg) 2\xi  \epsilon C^2 a_k^2
\end{multline*}
Summing from $m=0$ to $k$ and using Lemma \ref{lem2} gives:
\begin{multline*}
\sum_{m=0}^k 2a_k (\f(\vv_k) - f(w^*))\leq  d^2(\w_0 ,w^*)   
+ \Bigg\{  \frac{2\sqrt{n}MD}{1-\sigma_2(Q)} +  \xi  C^2    +  \frac{2}{ \Big( 1 - \frac{\mu_{\text{max} }\epsilon}{2} \Big)} \Bigg(\frac{(1 + \epsilon^2 L_{\varphi}^2)}{\epsilon }     + \frac{\mu_{\text{max}}}{2}\Bigg)2\epsilon \xi  C^2  \Bigg\} \sum_{m=0}^k a_k^2
\end{multline*}
\end{proof}
Using (\ref{f-bound}) and Lemma \ref{con}, we can get an upper bound on the objective values $f(v_k^i)$. We begin by noting that
$
 |f^i(v^i_k) -f(v^j_k) |\leq L d(v_k^i,v_k^j)  
$ for some constant $L$, so that 
\begin{multline}\label{f-boundll}
\sum_{m=0}^k 2a_k \big(f(v_m^i) - f(w^*)\big)\leq  2 \sum_{m=0}^k a_kL\sum_{j=1}^n  d(v_k^i,v_k^j)  
\\ + d^2(\w_0 ,w^*)    + \Bigg\{  \frac{2\sqrt{n}MD}{1-\sigma_2(Q)} +  \xi  C^2    +  \frac{2}{ \Big( 1 - \frac{\mu_{\text{max} }\epsilon}{2} \Big)}  \Bigg(\frac{(1 + \epsilon^2 L_{\varphi}^2)}{\epsilon }     + \frac{\mu_{\text{max}}}{2}\Bigg) 2\epsilon \xi  C^2  \Bigg\} \sum_{m=0}^k a_k^2
\end{multline}
which gives
\begin{multline*}
\inf_{0\leq m \leq k} f(v_m^i) \leq f(w^*) +L  \sup_{0\leq m \leq k} \sum_{j=1}^n  d(v_k^i,v_k^j) + \frac{1}{2\sum_{m=0}^k  a_m} \Bigg(  d^2(\w^i_0 ,w^*)   
+ \Bigg\{  \frac{2\sqrt{n}MD}{1-\sigma_2(Q)} +  \xi C^2   +\\ \frac{2}{ \Big( 1 - \frac{\mu_{\text{max} }\epsilon}{2} \Big)}   \Bigg(\frac{(1 + \epsilon^2 L_{\varphi}^2)}{\epsilon}     + \frac{\mu_{\text{max}}}{2}\Bigg) 2\epsilon \xi  C^2 \Bigg\} \sum_{m=0}^k a_k^2 \Bigg).
\end{multline*}
This implies, since $\sum_{k\geq 0} a_k^2 < \infty,\,\sum_{k\geq 0} a_k = \infty$ and $d(v_k^i,v_k^j) \to 0$, that
$$
\liminf_{k\to \infty} f(v_k^i) = f(w^*).\,\,\forall i,
$$

\subsection{Stochastic Gradient Optimization}
In this part, we analyze the algorithm for stochastic optimization. This means that we do not have access to the exact gradient but only an unbiased sample of it. An important motivating example of this scenario is the case when $f^i$ at each node is itself a sum of functions,
$$
f^i(x) = \frac{1}{n_i}\sum_{p =1}^{n_i} f^{p}(x)
$$
and to calculate the gradient, one may randomly uniformly sample a function $f^{\omega}$ from the set  $\{f^{p}\}_{p =1}^{n_i}$. This implies:
$$
\EX [\nabla f^{\omega} (x) | x] =\EX_{\omega} \, [\nabla f^{\omega}(x)]  = \ \frac{1}{n_i}\sum_{p=1}^{n_i} \nabla f^{p}(x),
$$
where $\EX[\cdot |x]$ denotes the condition expectation  given $x$. The above situation is quite common in machine learning with large distributed datasets so that stochastic optimization is preferable, or in online optimization where the samples are arriving in real time. Accordingly, in this section we assume that we have unbiased gradients with finite variance, i.e. for all $i$, at any instant $k$ we have access to an unbiased sample $\ff^{k,i}$ such that
\begin{equation}\label{stochastic}
\EXX\,  \nabla \ff^{k,i} (v^i_k) = \nabla f^i (v^i_k) \text{ and } \EXX \,\| \nabla \ff^{k,i} (\vv_k) - \nabla f^i(\vv_k)\|^2  \leq \sigma^2,
\end{equation} 
where $\EXX $ is the conditional expectation given the sigma algebra $\sigma(\vv_0,\cdots,\vv_k)$ generated by all the past information of the algorithm. We assume that $\{\ff^{i,k}\}_{k\geq 0}$ are conditionally independent without any loss of generality. E.g., for the case of sampling from a dataset, this would amount to sampling $f^\omega$ independently at each node for all $k$. We note that, the gradient tracking step (S2) in Algorithm 1 gets modified to:
\begin{equation}\label{alt}
\gd_{k}^i = \sum_{j=1}^{n}q_{ij} \T_{v_{k-1}^j}^{v_{k}^i}  \,\gd_{k-1}^j + \nabla {\tilde{f}}{}^{k,i}(v_{k}^{i}) - \T_{v_{k-1}^i}^{v_{k}^i} \nabla {\ff}{}^{k-1,i}(v_{k-1}^{i}) ,
\end{equation}
%We note that the above iteration entails the calculation of two gradients, $ \nabla {\tilde{f}}{}^{i,k}(v_{k}^{i})$ and $ \nabla {\ff}{}^{i,k}(v_{k-1}^{i})$, so it has marginally greater computational overhead than the previous case. We have the following result for the expectation and variance of $\gd_{k}^i$:
\begin{lem}\label{2-1}
For all $i$, we have the following relation
\begin{equation*}
\EXX \gd_{k} = g_k^i\,\, \text{  and  } \,\,\EX \|\gd^i_{k}-g_k^i\|^2 \leq \Big(\frac{3n+1}{1-q_{\text{max}}^2}\Big) \sigma^2,
\end{equation*}
where $q_{\text{max}} = \max_{i,j} q_{ij} <1$.
\end{lem}
\begin{proof}
Taking conditional expectations in (\ref{alt}), one can see from an easy induction argument that:
\begin{equation*}
\EXX \,\gd_{k}^i = \sum_{j=1}^{n}q_{ij} \T_{v_{k-1}^j}^{v_{k}^i}  \,g_{k-1}^j + \nabla {f}{}^{i}(v_{k}^{i}) - \T_{v_{k-1}^i}^{v_{k}^i} \nabla {f}{}^{i}(v_{k-1}^{i}) .
\end{equation*}
Let $\xi^{k,i} (v_k^i): =  \nabla  \ff^{k,i} (v_k^i) - \nabla f(v_k^i)$. We note using (\ref{stochastic}), that $\xi_{k,i}$ are zero mean and conditionally independent with the conditional variance, $\text{Var} (\xi^{k,i} (v_k^i)) \leq \sigma^2 $. Then, performing a recursion in (\ref{alt}) and subtracting from (\ref{11-00}) (see Appendix I), we have, 
\begin{align*}
 \gd_{k}^i - g_k^i &=  \sum_{j=1}^{n} q_{ij}^k\T_{v_{0}^j}^{v_{k}^i}  \xi^{0,i}  + \sum_{p=1}^{k-1} \sum_{j=1}^{n} q_{ij}^{k-p} \Big(\T_{v_{p}^j}^{v_{k}^i}   \xi^{p,j} - \T_{v_{p-1}^j}^{v_{p}^i} \xi^{p-1,j} \Big) + \xi^{k,i} - \T_{v_{k-1}^i}^{v_{k}^i} \xi^{k-1,i}\\ 
 &=  \sum_{p=0}^{k-2}\sum_{j=1}^{n} \Big(q_{ij}^{k-p} \T_{v_{p}^j}^{v_{k}^i}   \xi^{p,j} - q_{ij}^{k-p-1}\T_{v_{p}^j}^{v_{p}^i} \xi^{p,j} \Big) +  \sum_{j=1}^n\T_{v_{k-1}^j}^{v_{k}^i} q_{ij} \xi^{k-1,j}   +  \xi^{k,i} - \T_{v_{k-1}^i}^{v_{k}^i} \xi^{k-1,i}
\end{align*}
Recalling the fact that $\text{Var}(c_1 X_1 + c_2 X_2) = c^2_1\text{Var}(X_1)+ c^2_2 \text{Var}(X_2)$ for independent $X_1, \,X_2$, we have 
\begin{align*}
\EXX \| \gd_{k}^i - g_k^i\|^2 & \leq n  \sum_{p=1}^{k-1}(q_{ij}^{2k-2p}+q_{ij}^{2k-2p-2})\sigma^2  + (n+1)\sigma^2\\
& \leq \frac{2n\sigma^2}{1-q^2_{ij}} +(n+1) \sigma^2.
\end{align*}
which proves the result.
\end{proof}
We note that the previous lemma implies:
$$
\EX \|\gd_k^i \|^2 \leq R^2: = 2\Big(C^2 +\Big(\frac{3n+1}{1-q_{\text{max}}^2}\Big) \sigma^2\Big),
$$
since $\|g_k^i\|^2 \leq C^2$. We next prove the stochastic analog of Lemma \ref{lem1}:
\begin{lem}
We have the following  descent relation
\begin{align*}
\epsilon\,\EX\|  \text{grad} \, \varphi(\mathbf{w}_k)\|^2   &\leq \frac{1}{ \Big( 1 - \mu_{\text{max}} \epsilon \Big)}\big(\EX \varphi(\w_{k}) - \EX \varphi(\w_{k+1})\big)+   \frac{2}{ \Big( 1 - \frac{\mu_{\text{max} }\epsilon}{2} \Big)}  \Bigg(\frac{(1 + \epsilon^2 L_{\varphi}^2)}{\epsilon}     + \frac{\mu_{\text{max}}}{2}\Bigg)R^2 a_k^2
\end{align*}
\end{lem}
\begin{proof}
We note that
\begin{align*} 
 \varphi(\w_{k+1})  &\leq  \varphi(\vv_{k}) + a_k \langle  \text{grad}_{\vv_k} \varphi(\mathbf{v}_k),  \tilde{\g}_k \rangle   + \frac{\mu_{\text{max}}\| \tilde{\g}_k \|^2}{2} a_k^2.
\end{align*}  
Taking conditional expectations and using Lemma \ref{2-1}, we get 
\begin{align*} 
\EXX \varphi(\w_{k+1})  &\leq  \varphi(\vv_{k}) + a_k \langle  \text{grad}_{\vv_k} \varphi(\mathbf{v}_k),  \g_k \rangle   + \frac{\mu_{\text{max}}R^2 }{2} a_k^2.
\end{align*} 
Recalling the Taylor series expansion of $\varphi(\vv_k)$, we have
\begin{equation*}
\varphi(\textbf{v}_{k}) \leq   \varphi(\w_{k}) +  \|\text{grad}_{\textbf{w}_{k}} \varphi(\mathbf{w}_k) \|^2\,\epsilon\, \Big( 1 - \frac{\mu_{\text{max}} \epsilon}{2}\Big)
\end{equation*}
Thus, combining the previous two inequalities we get the following: 
\begin{equation*}
\EXX \varphi(\w_{k+1})  \leq \varphi(\w_{k}) -    \|\text{grad}_{\textbf{w}_{k}} \varphi(\mathbf{w}_k)) \|^2\,\epsilon\, \Big( 1 - \frac{\mu_{\text{max}} \epsilon}{2}\Big) + a_k \langle  \text{grad}_{\vv_k} \varphi(\mathbf{v}_k),  \g_k \rangle   + \frac{\mu_{\text{max}}R^2 }{2} a_k^2 
\end{equation*}
Taking full expectations in the above gives:
\begin{equation*}
\EX \varphi(\w_{k+1})  \leq \EX\varphi(\w_{k}) -   \EX \|\text{grad}_{\textbf{w}_{k}} \varphi(\mathbf{w}_k)) \|^2\,\epsilon\, \Big( 1 - \frac{\mu_{\text{max}} \epsilon}{2}\Big) + a_k\EX \langle  \text{grad}_{\vv_k} \varphi(\mathbf{v}_k),  \g_k \rangle   + \frac{\mu_{\text{max}}R^2 }{2} a_k^2 
\end{equation*}
The rest of the proof is exactly the same as that for Lemma \ref{lem1}.
\end{proof}
We are in a position to state the main result. The convergence analysis is the along the same lines as in the previous section with some important modifications. 
\begin{thm}
Assuming condition (\ref{stochastic}), we have with probability 1,
$$
\liminf_k f(v_k^i) = f(w^*),
$$
if $\epsilon \in (0,2\mu^{-1}_{\text{max}})$, $\sum_{k\geq 0}a_k = \infty$ and $\sum_{k=0}^\infty a_k^2 <\infty $.
\end{thm}
\begin{proof}
We begin by taking conditional expectations in the following bound:
\begin{align*}
d^2(w^i_{k+1},w^* ) &\leq d^2(v^i_k ,w^*) +2 a_k\langle \gd^i_k, \E^{-1}_{v_k^i} ( w^*) \rangle +
 \xi \| \E^{-1}_{ v^i_{k}} (w_{k+1}^i) \|^2 \\
\EXX d^2(w^i_{k+1},w^* )  &\leq d^2(v^i_k ,w^*) +2 a_k\langle g^i_k, \E^{-1}_{v_k^i} ( w^*) \rangle +\xi \| \E^{-1}_{ v^i_{k}} (w_{k+1}^i) \|^2 
\end{align*}
Taking full expectations, using the geodesic convexity property of  $f(\cdot)$, we have
\begin{align*}
\EX \,d^2(\w_{k+1}, \w^* )   &\leq \EX \,d^2(\vv^i_k ,\w^*) +2 a_k\Big(\EX (\f(\w_k ) -f(\w^*)) \Big) 
 + 2a_k n^{-1}D \sum_{i=1}^n \EX \| g^i_k - \nabla f(v_k^i)\| +  \xi  R^2 a_k^2 
\end{align*}
The analogous expression for $\vv_k$ gives
\begin{align*}
\EX d^2(\vv_{k},\w^* ) &\leq d^2(\w_k ,w^*) +  \xi  \epsilon^2 \EX \|  \text{grad}_{\textbf{w}_{k}} \varphi(\mathbf{w}_k) \|^2,
\end{align*}
where we use the convexity of $\varphi(\cdot)$ and the fact that $\varphi(\textbf{w}^*) = 0 $ with $\w^*:= (w^*,\cdots,w^*)$ and $\EX \varphi(\textbf{w}_k)\geq 0$.
Combining the previous two expressions, we have, after performing the same manipulations as Theorem 6:
\begin{multline*}
\sum_{m=0}^k 2a_k (\EX \f(v^i_m) - f(w^*))\leq   L_f \sum_{j=1}^n \EX d(v_k^i,v_k^j) + \frac{1}{2\sum_{m=0}^k  a_m} \Bigg( d^2(\w^i_0 ,w^*)   
+ \Bigg\{  \frac{2\sqrt{n}MD}{1-\sigma_2(Q)} +  \xi C^2   + \\  \frac{2}{ \Big( 1 - \frac{\mu_{\text{max} }\epsilon}{2} \Big)} \Bigg(\frac{(1 + \epsilon^2 L_{\varphi}^2)}{\mu_{\text{max}}}     + \frac{\mu_{\text{max}}}{2}\Bigg)2\epsilon \xi  R^2  \Bigg\} \sum_{m=0}^k a_k^2 \Bigg).
\end{multline*}
We have $\EX d(v_k^i,v_k^j) \to 0$, so  that 
$$
\liminf_k \EX \f(v^i_k) -f(w^*)  = 0.
$$
Since $f(v_k^i) \geq f(w^*)$, we can apply Fatou's Lemma to conclude
$$
0\leq \EX \liminf_k  \f(\vv_k) -f(w^*)\leq \liminf_k \EX \f(\vv_k) -f(w^*)  =0
$$
which implies $ \liminf_k  \f(\vv_k) -f(w^*)=0$ w.p.1.
\end{proof}
 \begin{figure*}
        \centering
        \begin{subfigure}[b]{0.490\textwidth}
            \centering
            \includegraphics[width=\textwidth]{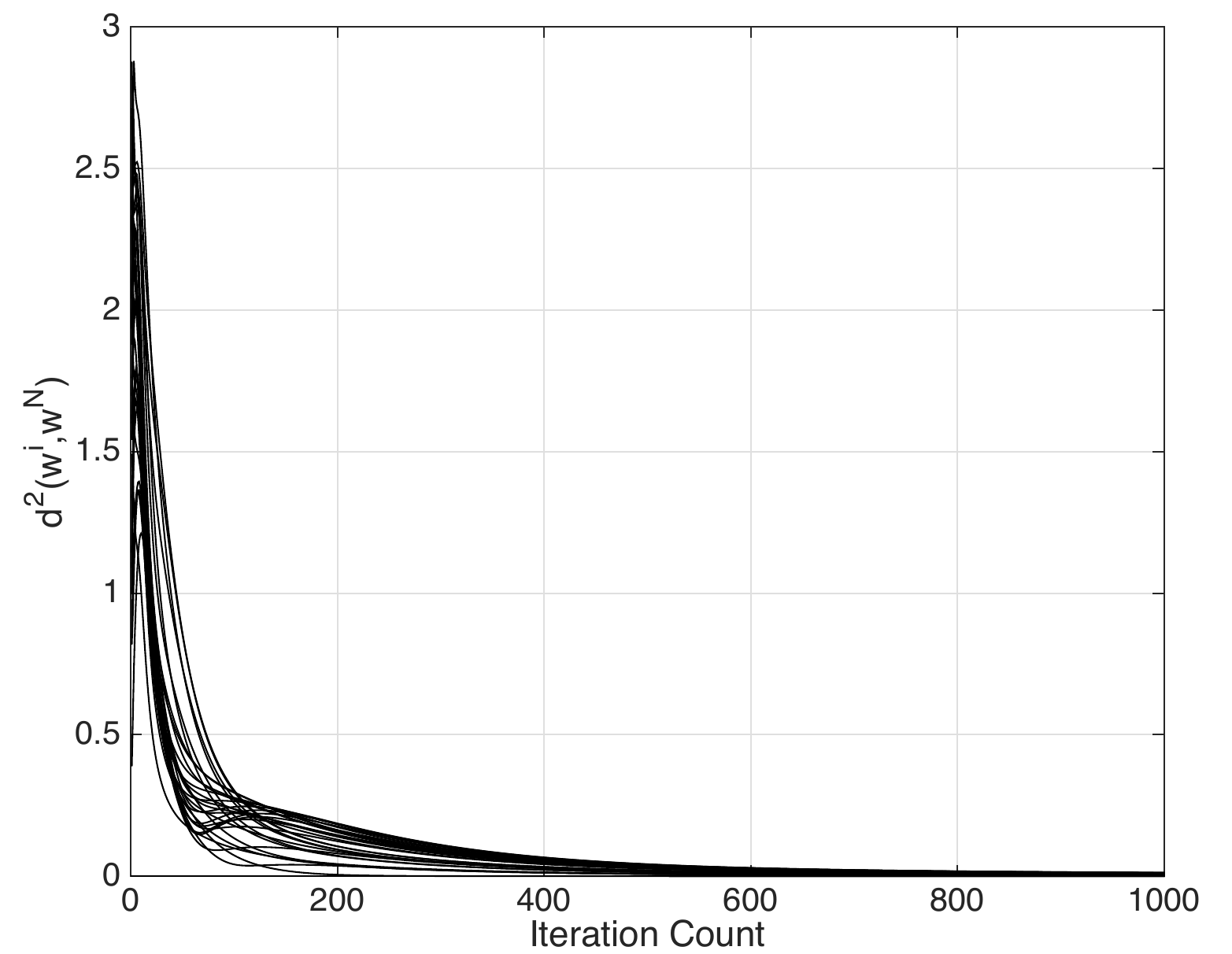}
            \caption[Network Graph]%
            {{\small Network Graph}}    
            \label{fig:mean and std of net14}
        \end{subfigure}
        \hfill
        \begin{subfigure}[b]{0.485\textwidth}  
            \centering 
            \includegraphics[width=\textwidth]{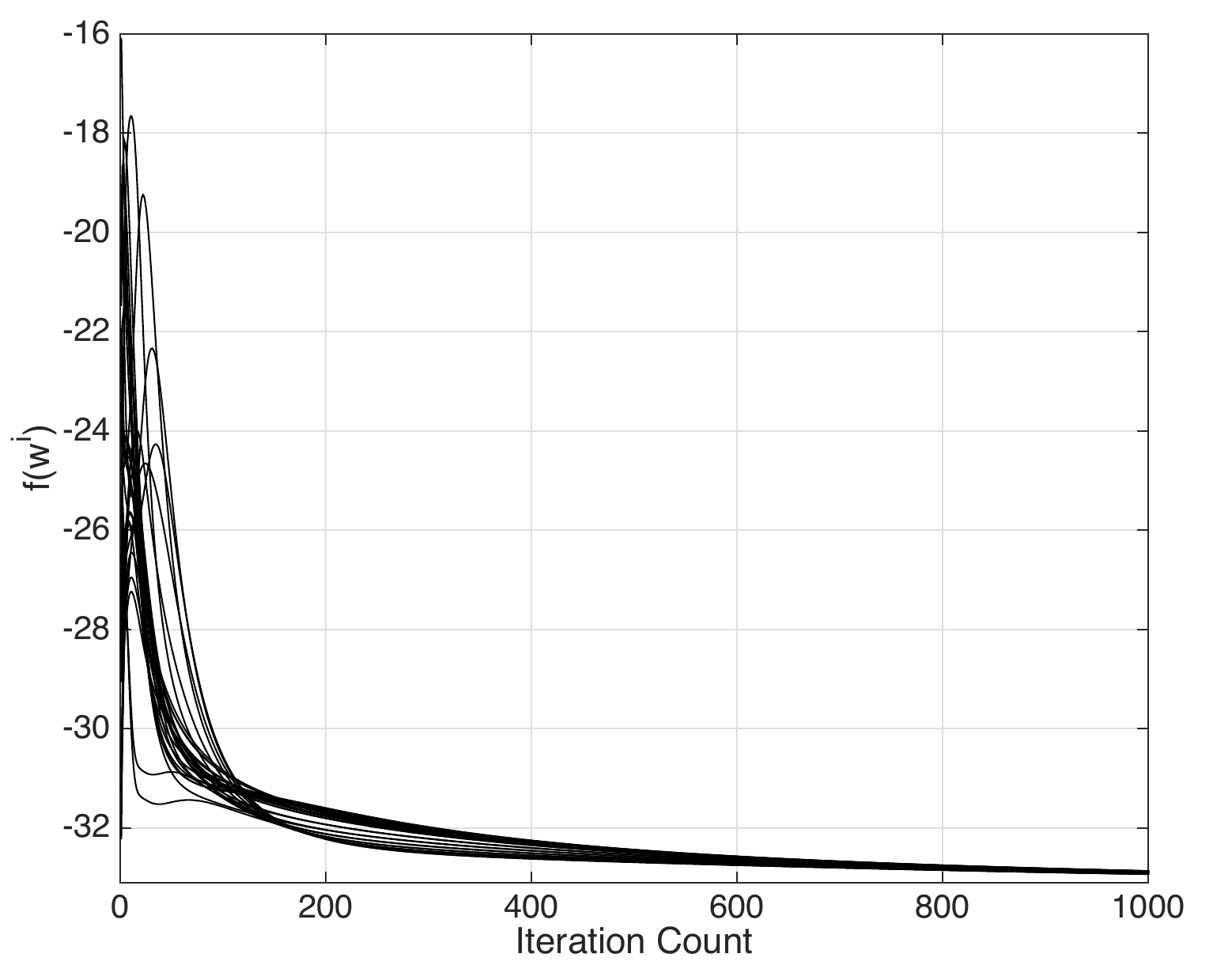}
            \caption[Distance between the agents vs Iteration count]%
            {{\small Distance between the agents vs Iteration count}}    
            \label{fig:mean and std of net24}
        \end{subfigure}
        \vskip\baselineskip
        \begin{subfigure}[b]{0.490\textwidth}   
            \centering 
            \includegraphics[width=\textwidth]{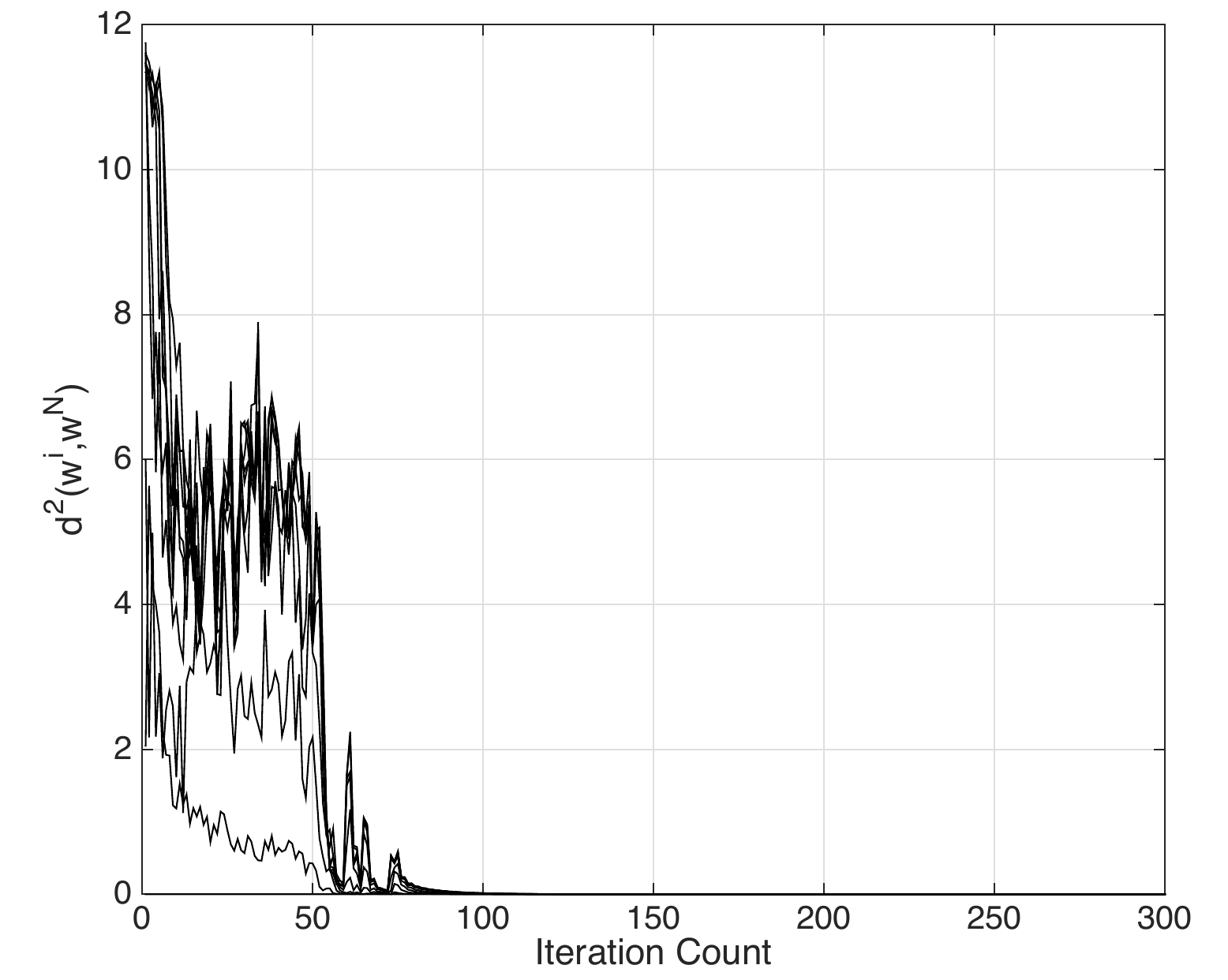}
            \caption[Function value vs Iteration Count]%
            {{\small Function value vs Iteration Count}}    
            \label{fig:mean and std of net34}
        \end{subfigure}
        \quad
        \begin{subfigure}[b]{0.47\textwidth}   
            \centering 
            \includegraphics[width=\textwidth]{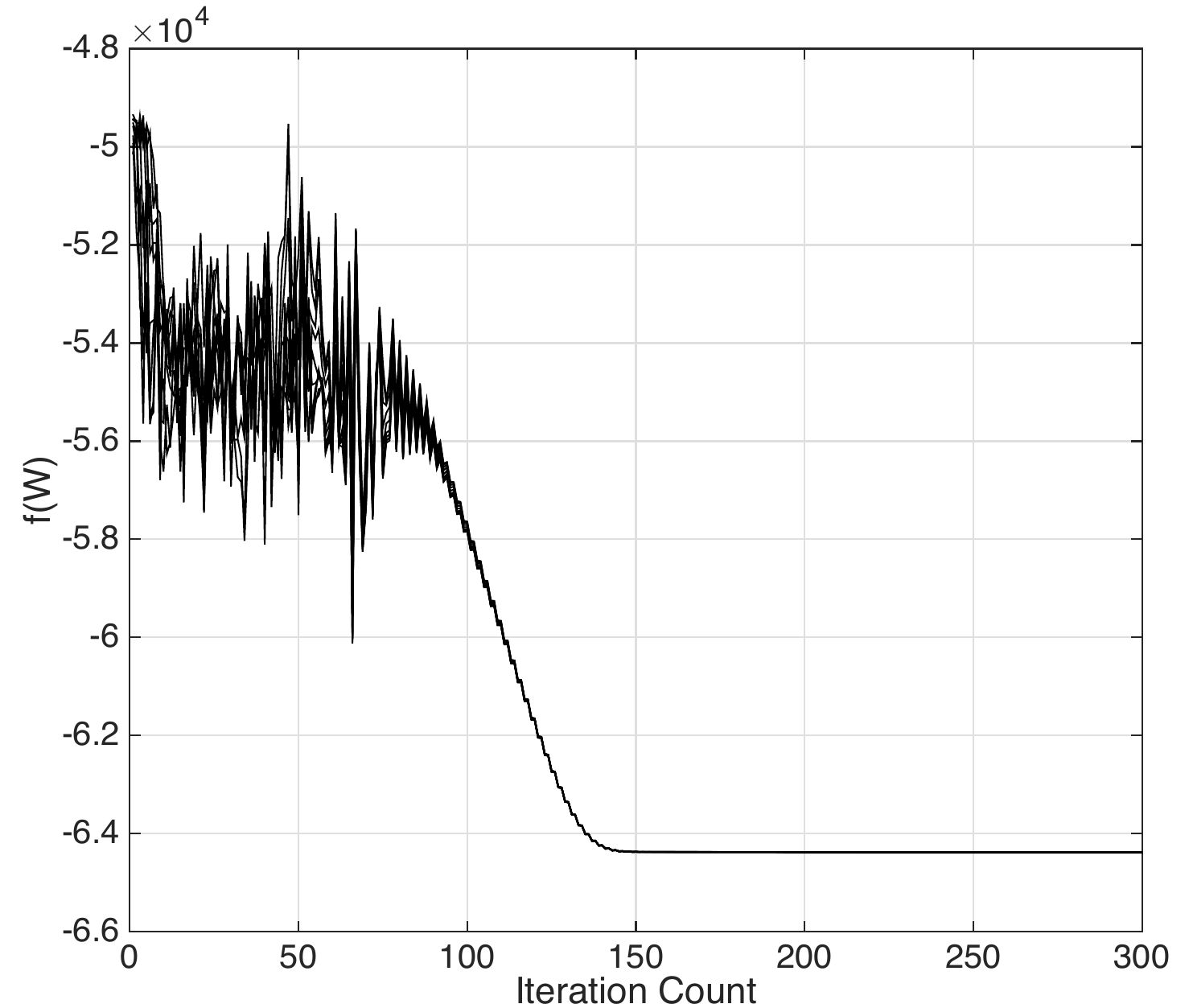}
            \caption[ Consensus Failure]%
            {{\small Consensus Failure}}    
            \label{fig:mean and std of net44}
        \end{subfigure}
        \caption[ The average and standard deviation of critical parameters ]
        {\small (a) This shows the plot of the distance between the various agents vs. the number of iterations. We fix an agent $i$ and plot the value $d^2(w^i,w^j)$ for all $j \in \mathcal{V}$. (b) This shows the function value, $f(w^i)$ for all $i$, against the iteration count. (c)-(d) The same results for Example 5.2. } 
        \label{fig:mean and std of nets}
    \end{figure*}
\section{ Numerical Experiments }
In this section we demontrate the algorithm on two important problems which can be solved within Riemannian framewok.  In particular, these examples demonstare the performence of proposed algorithm on the spherical and the Grassman manifolds. The algorithm is simulated in a MATLAB environment using the Manopt toolbox, \cite{manopt}.

\subsection{Computing the leading eigenvector}
The problem entails computing the leading eigenvector of data matrix whose entries are stored at different nodes of a communication network. The problem can be precisely stated as
\begin{align*}
\text{minimize} & \,\,\,\qquad \Big\{-w^t A w \,\,\overset{\Delta}{=} \,\, -w^T\Big(\sum_{i=1}^n z_iz_i^T \Big)w \Big\},\\
\text{subject to}& \,\,\,\qquad w\in \mathbb{S}^p, 
\end{align*} 
where $z_i \in \mathbb{R}^{p+1}$ is the data entry located at the \textit{i}'th node and $\mathbb{S}^p $ is the $n$-dimensional sphere defined as,
$$  \mathbb{S}^p=\{w\in\mathbb{R}^{p+1} \,:\, w^T w=1 \}.  $$ 
We note that the distributed version of the  problem could be  posed in the form (\ref{mainprob}) by letting $f^i(w) :=  -w^T z_iz_i^Tw$. The tangent plane $\T_w\SA^p$ at a point $w$ is given by,
$$ \T_w \mathbb{S}^p =\{ x \in \mathbb{R}^{n+1} \,:\,x^Tw =0\}.$$
The Riemannian metric is the usual inner product inherited from $\mathbb{R}^{p+1}$. The  sectional curvature is constant with $\kappa_{min}=\kappa_{max}=1$. Also, $\text{inj} \, \mathbb{S}^p = \pi$ and $ r_c = \pi/2$.    
The Riemmanian gradient $\text{grad }f(\cdot)$ for any function $f(\cdot)$ can be obtained by the orthogonal projection $\PA_{\T_w}(\cdot)$ of the Euclidean gradient $\nabla f(w)$ on the tangent space,
$$       \text{grad }_w\,f(w) =  \PA_{\T_w} ( \nabla f(w)) = (I_{p+1} - ww^T)\nabla f(w),       $$ 
where $I_{p+1}$ is the $(p+1) \times (p+1)$ identity matrix. The geodesic $t\to x(t)$ expressed as a function of $x_0 \in S^{p-1}$ and $\dot{x}_0 \in \T_{x_0}S^{p-1}$ is given by (Example 5.4.1, \cite{Absil}) :
$$x(t)=x_0\text{cos}(\|\dot{x}_0t\|) +\dot{x}_0\frac{1}{\| \dot{x}_0\|}\text{sin}(\|\dot{x}_0\|t).$$
An alternative to this is provided by approximating the geodesic with a retraction which just involves normalizing (i.e., divide by the norm) $w_n$ at each step. The Riemannian distance between any two points is given as $$d(x,y) = \cos^{-1}(x^Ty).$$ The vector transport associated along the geodesic is given by the following expression:
\begin{equation*}     
\T_{x_0}^{x(t)}(\xi_x) = \Big\{  u  \cos(\|\dot{x}_0\|t) x^T(0)  - \sin(\|\dot{x}_0\|t)xu^T  +(I_{p+1}-uu^T) \Big\} \xi_x,       
\end{equation*}  
where  $u =\frac{\dot{x}_0}{\|\dot{x}(0)\|}$.

The matrix $Q$ here is generated using Metropolis weights and the number of nodes is set to $n=30$. We let $p=3$. To generate the initial conditions, a random point $w_0 \in\mathbb{S}^3$ is selected and the initial conditions at each node are set to $w^i_0 = \E_{w_0}(v^{i})$, where $\{v^i\}_{i=1}^N$ are tangent vectors drawn from an isotropic Gaussian distribution with standard deviation $\sigma = 0.2$. The algorithm is run for $k=10^3$ iterations (although convergence is achieved pretty quickly due to a small matrix size) with step seizes $\epsilon = 0.1$ and $ \,a_k= 1/k$. The results are plotted in  Figure 1.The agents achieve consensus relatively early and stay in that configuration for the rest of the run. We remark that the number of iterations required typically increases with the dimension of the problem. Figure 1b shows the function value, $f(w^i)$ for all $i$, against the iteration count. Note that as the agents achieve consensus, the values decrease in cohesion for all agents.
 
\subsection{Principal Component Analysis }

We next consider a generalization of the previous problem known as Principal Component Analysis (PCA). PCA entails the computation of the $r$ principal eigenvectors of a $p \times p$ covariance
matrix $A$ given by $$A = \mathbb{E}[z_{k}z_{k}^{T}]$$ 
with $z_{1},...,z_{k},...$ being a stream of uniformly bounded $p$-dimensional data vectors. The cost function for PCA is :
\[
C(W)=-\frac{1}{2}\mathbb{E}[z^{T}W^{T}Wz]=-\frac{1}{2}\text{Tr}(W^{T}AW),
\]
where $\text{Tr}(\cdot)$ denotes the trace of matrix and $W$ belongs to the Stiefel manifold ${S}_{p,r}$ defined as : 

\[
\mathcal{S}_{p.r}=\{X\in\mathbb{R}^{p\times r}\,:\,X^{T}X=I_{r}\}.
\]
 The cost function is invariant to the transformation $W \mapsto WO$ where $O\in \mathcal{O}(r)$ which represents the orthogonal group. Thus the state space can identified with the Grassmann manifold :
$$
\mathcal{G}(p,r) = \big\{ \mathcal{S}_{p,r} / O(r) \big\}, 
$$
so that any $[G] \in \mathcal{G}(p,r)$ is the equivalence class, $$ [G]=\{ WO,\,O\in\mathcal{O}_r \}.$$   
For details on the geometry of this manifold we refer the reader to \cite{Edel}. The Riemannian gradient of $C(W)$ under the sample $z$
is given by, $$ \text{grad}_W \,f(W)=(I_n-WW^{T})zz^{T}W.$$ A retraction  which could be used here is given by (Example 4.1.3, \cite{Absil}) :
$$ \mathcal{R}_{W}(a H) = \text{qf}(W + aH),$$
where $\text{qf}(\cdot)$ gives the orthogonal  factor in the QR-decomposition of its argument. The formula for the parallel transport on Grassmann manifolds is given in (Example 8.1.3 \cite{Absil}) and for vector transport in (Example 8.1.10 \cite{Absil}).

We consider a synthetic dataset of $d=10^4$ measurements with dimension of each measurement being $p=10^3$. The measurements are assumed to be distributed throughout a network of $10$ nodes. The connectivity graph is generated in a similar manner as the previous example and the results are plotted in Figure 1(c)-(d). A point to note here is that consensus, as in the previous case, is achieved pretty early around $k=100$ while the function value evaluated at the estimates keeps on decreasing till $k=150$.

%
%%
%\begin{figure}[H]
%\begin{center}
%\includegraphics[width=11cm,height=8cm]{fig1}
%\end{center}
%
%\end{figure}
%
%
%\begin{figure}[H]
%\begin{center}
%\includegraphics[width=10.5cm,height=8.5cm]{fig2}
%\end{center}

%\begin{figure*}[h]
%\centering
%    \includegraphics[width=0.47\linewidth]{fig1}\hfil
%    \includegraphics[width=0.47\linewidth]{fig2}
%  \caption{caption here}
%\end{figure*}

\section*{Appendix I: Proof of Lemma 3}\label{A2}
We recall the iteration for computing $g_k^i$:
\begin{equation}\label{g-track}
g_{k}^i =\sum_{j=1}^{n}q_{ij} \T_{v_{k-1}^j}^{v_{k}^i}  \,g_{k-1}^j + \nabla {f}{}^{i}(v_{k}^{i}) - \T_{v_{k-1}^i}^{v_{k}^i} \nabla {f}{}^{i}(v_{k-1}^{i}) ,
\end{equation}
where  $g_0^i = \nabla f^{i}(v_{0}^{i})$. Let $q_{ij}^m$ denote the $ij$'th entry of the matrix product $Q^m:= Q\times\cdots \times Q$. Doing a recursion on the expression for $g_k^i$ gives
\begin{multline}\label{11-00}
g_{k}^i =\sum_{j=1}^{n} q_{ij}^k\T_{v_{0}^j}^{v_{k}^i}   g_{0}^i +    \sum_{p=1}^{k-1} \sum_{j=1}^{n} q_{ij}^{k-p} \Big(\T_{v_{p}^j}^{v_{k}^i}  \big( \nabla {f}{}^{j}(v_{p}^{j}) - \T_{v_{p-1}^j}^{v_{p}^i} \nabla {f}{}^{j}(v_{p-1}^{j})\big) \Big) + \nabla {f}{}^{i}(v_{k}^{i}) - \T_{v_{k-1}^i}^{v_{k}^i} \nabla {f}{}^{i}(v_{k-1}^{i}),
\end{multline}
 where we have used the associative property of parallel transport, namely $\T^{x}_y \T^y_z = \T^y_z$. We show that the auxiliary sequence $\{g^i_k\}$ tracks the average gradient closely enough at each node. 
 
Next, we let $\bar{g}_k^i = \frac{1}{n}\sum_{j=1}^{n} \T_{v_k^j}^{v_k^i}  \,g_{k}^j $. We can write the expression for $\bar{g}_k^i$ using the definition of $\,g_{k}^j $ as 
\begin{align*}
\bar{g}_{k}^i &= \frac{1}{n}\sum_{j=1}^n \T_{v_{k}^j}^{v_k^i} \Big(\sum_{p=1}^{n} q_{jp} \T_{v_{k-1}^p}^{v_{k}^j}  \,g_{k-1}^p + \nabla {f}{}^{j}(v_{k}^{j}) - \T_{v_{k-1}^j}^{v_{k}^j} \nabla {f}{}^{j}(v_{k-1}^{j}) \Big) \nonumber\\
&= \frac{1}{n}\sum_{j=1}^n \sum_{p=1}^{n} q_{jp} \T_{v_{k-1}^p}^{v_{k}^i}  \,g_{k-1}^p + \frac{1}{n}\sum_{j=1}^n  \T_{v_{k}^j}^{v_{k}^i} \Big( \nabla {f}{}^{j}(v_{k}^{j})  - \T_{v_{k-1}^j}^{v_{k}^j} \nabla {f}{}^{j}(v_{k-1}^{j}) \Big), \nonumber \\
& = \frac{1}{n}\sum_{j=1}^n  \T_{v_{k-1}^p}^{v_{k}^i}  \,g_{k-1}^j + \frac{1}{n}\sum_{j=1}^n  \T_{v_{k}^j}^{v_{k}^i} \Big( \nabla {f}{}^{j}(v_{k}^{j}) - \T_{v_{k-1}^j}^{v_{k}^j} \nabla {f}{}^{j}(v_{k-1}^{j}) \Big), \nonumber \\
& = \T_{v_{k-1}^i}^{v_{k}^i}  \,\bar{g}_{k-1}^i + \frac{1}{n}\sum_{j=1}^n  \T_{v_{k}^j}^{v_{k}^i} \Big( \nabla {f}{}^{j}(v_{k}^{j}) - \T_{v_{k-1}^j}^{v_{k}^j} \nabla {f}{}^{j}(v_{k-1}^{j}) \Big),
\end{align*}
where we use the fact that $\sum_{j=1}^n q_{ij}=\sum_{i=1}^n q_{ij}=1$ in the third equality. Since $\bar{g}_{0}^i  = \frac{1}{n}\sum_{p=1}^n  \nabla {f}{}^{p}(v_{0}^{p}) $, using the previous equation, one can see by induction that $\bar{g}_{p}^i  :=\frac{1}{n}\sum_{p=1}^n   \nabla {f}{}^{p}(v_{p}^{p}) $ for any $p\geq 1$. Doing a recursion for $\bar{g}_{k-1}^i$ in the previous equation we have
\begin{equation}\label{2200}
\bar{g}_{k}^i  = \T_{v_{0}^i}^{v_{k}^i}  \,\bar{g}_{0}^i + \frac{1}{n} \sum_{p=1}^k\sum_{j=1}^n  \T_{v_{p}^j}^{v_{k}^i} \Big( \nabla {f}{}^{p}(v_{p}^{j}) - \T_{v_{p-1}^j}^{v_{k}^j} \nabla {f}{}^{p}(v_{p-1}^{j}) \Big).
\end{equation}
Subtracting  (\ref{2200}) from (\ref{11-00}), we get 
\begin{multline}\label{t-t-t}
 g_{k}^i - \bar{g}_{k}^i =\sum_{j=1}^{n} q_{ij}^k\T_{v_{0}^j}^{v_{k}^i}   g_{0}^i  -\T_{v_{0}^i}^{v_{k}^i}  \,\bar{g}_{0}^i  
 \sum_{p=1}^{k-1}\sum_{j=1}^n \Big( q_{ij}^{k-p} - \frac{1}{n} \Big) \T_{v_{p}^j}^{v_{k}^i} \Big( \nabla {f}{}^{j}(v_{p}^{j}) - \T_{v_{p-1}^j}^{v_{k}^j} \nabla {f}{}^{j}(v_{p-1}^{j}) \Big) \\+\ \nabla {f}{}^{i}(v_{k}^{i}) - \T_{v_{k-1}^i}^{v_{k}^i} \nabla {f}{}^{i}(v_{k-1}^{i})  - 
  \frac{1}{n}\sum_{j=1}^n  \T_{v_{k}^j}^{v_{k}^i} ( \nabla {f}{}^{j}(v_{k}^{j}) - \T_{v_{k-1}^j}^{v_{k}^j} \nabla {f}{}^{j}(v_{k-1}^{j}) ).
\end{multline}
From the mixing bound (\ref{mix-bound}), we have
\begin{multline}\label{tetrasumi}
 \|g_{k}^i - \bar{g}_{k}^i \|= \sigma^k_2(Q)C \sqrt{n} + \sum_{p=1}^{k-1}  \sigma^{k-p}  \| \nabla {f}{}^{j}(v_{p}^{j})- \T_{v_{p-1}^j}^{v_{k}^j} \nabla {f}{}^{j}(v_{p-1}^{j}) \| +  \| \nabla {f}{}^{i}(v_{k}^{i}) - \T_{v_{k-1}^i}^{v_{k}^i} \nabla {f}{}^{i}(v_{k-1}^{i}) \\ - 
   \frac{1}{n}\sum_{j=1}^n  \T_{v_{k}^j}^{v_{k}^i} \Big( \nabla {f}{}^{j}(v_{k}^{j}) - \T_{v_{k-1}^j}^{v_{k}^j} \nabla {f}{}^{j}(v_{k-1}^{j}) \Big)\|.
\end{multline}
We next consider bounding $ \| \nabla {f}{}^{j}(v_{p}^{j}) - \T_{v_{k-1}^j}^{v_{p}^j} \nabla {f}{}^{j}(v_{p-1}^{j})\|$ for $1\leq p\leq k$. We begin by noting that for any $i=1,\cdots,n$, we have
\begin{equation}\label{a-kbound}
\|\nabla {f}{}^{i}(v_{k}^{i}) - \T_{v_{k-1}^i}^{v_{k}^i} \nabla {f}{}^{i}(v_{k-1}^{i}) \| \leq d(v^i_k,v^i_{k-1}) 
\end{equation}
 By triangle inequality we have,
\begin{align*}
d(v^i_k,v^i_{k-1}) &\leq d(v^i_k,w^i_{k}) + d(w^i_k,v^i_{k-1}) \\
&  \leq \epsilon\| \text{grad}_{w_k^i} \varphi(\w_k)\| + a_{k-1} \| g^i_{k-1}\|,
\end{align*}
Thus, we have 
\begin{align*}
&a_k\|\nabla {f}{}^{i}(v_{k}^{i}) - \T_{v_{k-1}^i}^{v_{k}^i} \nabla {f}{}^{i}(v_{k-1}^{i}) \|  \leq a_k d(v^i_k,v^i_{k-1})\\
&\leq a_k \big( \epsilon\| \text{grad}_{w_k^i} \varphi(\w_k)\| + a_{k-1} \| g^i_{k-1}\|\big)
\leq  \Big( \frac{a_k^2}{2} + \epsilon^2 \frac{\| \text{grad}_{w_k^i} \varphi(\w_k) \|^2}{2}+ aC a_k^2\Big),
\end{align*}
where we have used $a_{k-1} \leq a\, a_k$. Using the above in (\ref{t-t-t}), 
\begin{align}\label{1.1}
 a_k\|g_{k}^i - \bar{g}_{k}^i \|&\leq  \sqrt{n} \sum_{m=0}^{k-1}  \sigma_2^{k-m} (Q)\Big( \frac{a_m^2}{2} + \epsilon^2 \frac{\| \text{grad}_{w_m^i} \varphi(\w_m) \|^2}{2}   + aC a_m^2 \Big)  +\\
 &\qquad \qquad \qquad  \Big(1+\frac{1}{n}\Big) \Big(\frac{a_k^2}{2} + \epsilon^2 \frac{\| \text{grad}_{w_k^i} \varphi(\w_k) \|^2}{2}  + aC a_k^2\Big ) \nonumber\\
&\leq  \sqrt{n} \sum_{m=1}^{k}  \sigma^{k-m}_2(Q) \Big( \frac{a_m^2}{2} + \epsilon^2 \frac{\| \text{grad}_{w_m^i} \varphi(\w_m) \|^2}{2}+ aC a_k^2\Big), 
\end{align}
where we assume, without loss of generality, that $1 \leq  a \,a_0^2$ to absorb the first term in (\ref{tetrasumi}) into the summation and $\sqrt{n} \geq  1+n^{-1}$ to absorb the last term in the above inequality. We also have
\begin{align}\label{2.1}
a_k\|g_{k}^i - \nabla f (v^i_k) \| \leq a_k\|g_k^i  - \bar{g}_k^i\| + a_k\|\bar{g}_k^i - \nabla f (v^i_k) \|
\end{align} 
The second term in the previous equation can be bounded as
\begin{align*}
 a_k\|\bar{g}_k^i - \nabla f (v^i_k) \| &\leq \frac{a_k^2}{2} +\frac{\|\bar{g}_k^i - \nabla f (v^i_k) \|^2}{2}\\
 &= \frac{a_k^2}{2} +\frac{\| \frac{1}{n}\sum_{p=1}^n   \nabla {f}{}^{p}(v_{k}^{j})   - \nabla f (v^i_k) \|^2}{2}\\
 &\leq \frac{ a_k^2}{2}   + \frac{1}{2n} \sum_{j=1}^n d^2(v^i_k,v^j_k)\\
 &\leq \frac{ a_k^2}{2}   + \frac{3}{2n} \sum_{j=1}^n (d^2(v^i_k,w^i_k) +  d^2(w^i_k,w^j_k) +d^2(w_k^j,v^j_k))\\
 &\leq \frac{ a_k^2}{2}   + \frac{9}{2} \| \text{grad} \, \varphi(\w_k)  \|^2, 
\end{align*} 
where we have used the fact that $d^2(w^i_k,w^j_k) \leq \| \text{grad} \, \varphi(\w_k)  \|^2$ by definition.  Using the above and (\ref{1.1}) in (\ref{2.1}), 
\begin{equation*}
 a_k\|g_{k}^i - \bar{g}_{k}^i \|=  \sqrt{n} \sum_{m=0}^{k}  \sigma_2^{k-m}(Q) \Big( a_m^2 + \epsilon^2\Big(1+\frac{9}{2\epsilon^2} \Big) \| \text{grad}\,\varphi(\w_m) \|^2  + aC a_m^2 \Big) 
\end{equation*}
and using Lemma 2 to substitute for $\| \text{grad} \, \varphi(\w_m) \|^2$ gives
\begin{multline*}
 a_k\|g_{k}^i - \bar{g}_{k}^i \|=  \sqrt{n}\Bigg\{1 +aC + \frac{2\epsilon\Big(1+\frac{9}{2\epsilon^2} \Big)}{1-\frac{\mu_{\text{max}}\epsilon}{2} } \Big( \frac{(1 + \epsilon^2 L_{\varphi}^2)}{ \epsilon}     + \frac{\mu_{\text{max}}}{2}\Big)C^2 \Bigg\}  \times\\ \sum_{m=0}^{k} \big( \sigma_2^{k-m}(Q) a_m^2 \big) + \frac{2\sqrt{n}\epsilon\Big(1+\frac{9}{2\epsilon^2}\Big)}{(1-\frac{\mu_{\text{max}}\epsilon}{2} )}\sum_{m=0}^k( \varphi(\w_m) - \varphi(\w_{m-1}))
\end{multline*}
Set $M$ equal to the constant inside the curly brackets in the above equation. This gives
\begin{equation*}
 \sum_{p=0}^k a_p\|g_{p}^i - \bar{g}_{p}^i \|=  \sqrt{n}M \sum_{p=0}^k \sum_{m=0}^{p}  \sigma_2^{p-m}(Q) a_m^2 
\end{equation*}
We note that the summation on the right can be written as:
$$
 \sum_{p=0}^k \sum_{m=0}^{p}  \sigma^{p-m} a_m^2  \leq \sum_{m=0}^k a^2_m \sum_{p=0}^{k-m} \sigma^p \leq \frac{1}{1-\sigma}  \sum_{m=0}^k a^2_m.
$$
Plugging the above bound in the previous inequality proves the result.

\section*{ Appendix II: Proof of Lemma 4}
The proof uses techniques similar to Theorem 27, \cite{Tron}. We prove the result for a constant sectional curvature $\Delta$. The result for a varying (but bounded) curvature can be established along similar lines by using the Rauch comparision theorems (see e.g. Theorem IX.2.2, IX.2.3, \cite{chavel}). We consider the function $\varphi_{i,j}: \M \times\M \to \mathbb{R}$ defined by $\varphi_{i,j}(x^i, x^j ) :=\frac{1}{2} d^2(x^i,x^j)$. Since $\varphi(\x) = \sum_{(i,j)\in \mathcal{E}} \varphi(x^i,x^j)$, it is enough we prove the result for $\varphi_{i,j}$. We drop the indices $i,j$ in the proof to avoid notational clutter.

Let $(x,\,y) \in \M^2$. Define the geodesic $\gamma (-\epsilon, \epsilon) \to \M^2$ with $\gamma(0) = x$ and $\gamma(1)=y $. The second order Taylor expansion (see e.g., eq. (5.26), Section 5.3, \cite{boumal}) of $\varphi$ gives:
\begin{equation} \label{taylor}
\varphi(y) = \varphi(x) + \langle \text{grad} \varphi(x), \E^{-1}_x(y) \rangle + \frac{1}{2} \langle \text{Hess} \varphi(\gamma(t))[  \dot{\gamma}_t],\dot{ \gamma}(t) \rangle  + \frac{1}{2} \langle  \text{grad} (\varphi(\gamma(t))), \ddot{\gamma}(t)\rangle 
\end{equation}
The last term can be ignored since the acceleration of a geodesic is zero. We show that the third term is positive for any point $ \gamma(t) =(\xx,\yy$) and $\dot{ \gamma}(t)  = (v_1,v_2)$ with $t \in (0,1)$.

Accordingly, consider any point $(\xx,\yy) \in \M^2$ with $d(\xx,\yy)\leq r <r^*$. Define the geodesics $\gamma_i :(-\epsilon,\epsilon) \to \M$, with $\gamma_1(0)=\xx,\,\dot{\gamma}_1(0) = v_1$ and $\gamma_2(0)=\yy,\,\dot{\gamma}_2(0) = v_2$. Let $\gamma_{12,t}(s)$ denote the geodesic joining the point $\gamma_1(t)$ to $\gamma_2(t)$. We note that since we restrict ourselves to radii less than the convexity radius, such a construction is possible without encountering any conjugate point pairs on $\gamma_{12}$. Also, let $\phi(t):=L(\gamma_{12,t}) $, where $L$ denotes the length of the geodesic segment, so that $ \phi(t) :=d (\gamma_1(t),\gamma_2(t))$ and $\phi(0)=r$. Thus, we need to calculate the second order derivative of $\phi^2(t)$ at $t_0=0$.

 To proceed, we define the geodesic variation as $\alpha : [0,1]\times [a,b] \to \M $, so that the map $s \to \alpha(t_0,s)$ traces the (normalized) geodesic $\gamma_{1,2,t_0}(s)$ with $\| \dot{\gamma}_{1,2,t}\|=1$ . Define the vector field $X(t):= \frac{\partial \alpha}{\partial t}:=D\alpha \frac{\partial}{\partial t}|_{t=t_0}$. We note that  $X(t)$ is a Jacobi field (page 36, \cite{sakai}) along $\gamma_{12,t}$. The double derivative of $\phi^2(\cdot)$ can be calculated by a  straightforward application of the second variation of arc length formula (Theorem II.4.3, \cite{chavel}) and is given by (Theorem 27, \cite{Tron}):
 \begin{equation}\label{dder}
\frac{d}{dt^2} \phi^2(t)\Big|_{t = t_0} =  r \langle  X^\perp (s), X^{\perp}(s) \rangle\big|_{a}^b
+ \Big( \Big\langle X(s), \frac{\gamma_{12,t_0}(s)}{\| \gamma_{12,t_0}(s)\|} \Big\rangle \Big|_{a}^b \Big)^2, 
 \end{equation}
where $X^{\perp}(t):= X(t) - \langle X(t),  \dot{\gamma}_{12,t} \rangle \dot{\gamma}_{12,t}$ is the part of $X$ orthogonal to $\dot{\gamma}_{12,t}$. We consider the first term and show that it is non-negative to establish the result. Since the initial condition $X(0)$ is not zero, to use the standard Jacobi field theory we split $X(t)$ as $X(t) =X_1(t) + X_2(t)$, where 
\begin{equation}\label{i.c.}
X_1(0)=0 \,; X_2(0) =X^{\perp}(0);  \,X_1(r)=X^{\perp}(r); \,X_2(r) =0.
\end{equation}
Such a (unique) decomposition is possible (see e.g. Lemma 2.4, \cite{sakai}). Then, the solution to the Jacobi equation with the  above initial conditions, for a manifold with constant curvature, is given as (see Section 2.3, Chapter 5, \cite{dcarmo} or Section 3, \cite{chavel}):
\begin{equation}\label{h.def}
X_i(t) =  S_{\Delta}(t) E_i(t)  \text{ and } \nabla X_i(t) =  C_{\Delta}(t) E_i(t) , \,\,i=1,\,2
\end{equation}
where $E_i(t)$ is some field parallel to the $\gamma_{12,t}$ with $\| E_i(t)\| =1 $ and 
\begin{multline*}
 C_{\Delta}(t):=  \cos (\sqrt{\Delta}t),\,\, S_{\Delta}(t) :=  \frac{1}{\sqrt{\Delta}}\sin (\sqrt{\Delta}t),\text{ if } \Delta >0; \qquad 
  C_{\Delta}:= t,\,\,  S_{\Delta}= 1,\text{ if } \Delta =0;\\
    C_{\Delta}(t):=  \cosh (\sqrt{|\Delta|}t),\,\,\, S_{\Delta}(t) :=  \frac{1}{\sqrt{|\Delta|}}\sinh (\sqrt{|\Delta|}t),\text{ if } \Delta < 0.  
\end{multline*}
We also have the easily verifiable property (see e.g. Theorem 22, \cite{Tron}),
\begin{equation}\label{h.1}
 \frac{\langle \nabla X(t),X(t)\rangle }{\langle X(t) ,X(t)\rangle } = \frac{C_{\Delta}  (t)}{S_{\Delta}(t)} \text{ and } \| \nabla X (0)\| \leq \frac{\|X(t)\|}{S_\Delta(t)}
 \end{equation}
 We have from  (\ref{dder}):
\begin{align*}
\frac{1}{r}\frac{d}{dt^2} \frac{\phi^2(t)}{2}  &\geq \langle \nabla X_1(r),  X_1(r) \rangle + \langle \nabla X_2(r), X_1(r) \rangle  - \langle  \nabla X_1(0), X_2(0) \rangle  -\langle \nabla X_2(0), X_2(0) \rangle
\end{align*}
 We remark that since the initial conditions for $X_2(t)$ are reversed (i.e. $X_2(r)=0$ and $X_2(0)=X^{\perp}(0)$), when considering $\nabla X_r(0)$, we have to parametrize $\gamma_{12,t_0}(s)$ as $s' =r-s$. This gives $\nabla X(r)=-\nabla X(0)$. We first consider the case $\Delta \leq 0$. We have from using  (\ref{h.1}) in the the previous equation (alongwith $-\nabla X_2(0)=\nabla X_2(r)$):
 \begin{align*}
\frac{1}{r}\frac{d}{dt^2} \frac{\phi^2(t)}{2}  & \geq  \frac{C_\Delta(r)}{S_\Delta(r)} \|X_1(r) \|^2 - \frac{1}{S_\Delta}\|X_2(0)\| \|X_1(r)\| - \frac{1}{S_\Delta(r)} \|X_1(r)\| \| X_2(0)\|  + \frac{C_\Delta(r)}{S_\Delta(r)}  \|X_2(0) \|^2 \\
& \geq  \frac{C_\Delta(r)}{S_\Delta(r)} \|X_1(r) \|^2 + \frac{C_\Delta(r)}{S_\Delta(r)}  \|X_2(0) \|^2  - \frac{1}{S_\Delta} \Big(  \|X_1(r) \|^2 +  \|X_2(0) \|^2 \Big)\\
& =  \frac{1}{S_\Delta(r)}\Bigg(  C_\Delta(r)- 1 \Bigg) \Big( \|X^{\perp}(0)\|^2 +  \|X^{\perp}(r)\|^2 \Big),
\end{align*}
after using initial conditions (\ref{i.c.}).The right hand side can be seen to be positive for all $r>0$ using the definition of $C_\Delta$. 
 
 We next consider $\Delta > 0$. Set $a:= \|X^{\perp}(0)\|$ and $b=\|X^{\perp}(r)\|$. We have:
 \begin{align*}
\frac{1}{r}\frac{d}{dt^2} \frac{\phi^2(t)}{2}  & \geq  \frac{C_\Delta(r)}{S_\Delta(r)} \|X_1(r) \|^2 - C_\Delta(0)\|E_2(0)\| \|X_1(r)\| -C_\Delta(0)\|E_1(0)\| \| X_2(0)\|  + \frac{C_\Delta(r)}{S_\Delta(r)}  \|X_2(0) \|^2 \\
& \geq  \frac{C_\Delta(r)}{S_\Delta(r)} \|X^{\perp}(r) \|^2 - \|X^{\perp}(r)\| -\| X^{\perp}(0)\|  + \frac{C_\Delta(r)}{S_\Delta(r)}  \|X^{\perp}(0) \|^2 \\
&=  (a^2+b^2) \frac{C_\Delta(r)}{S_\Delta(r)}  \Bigg( 1- \frac{ \tan (\sqrt{\Delta}r) (a+b)}{\sqrt{\Delta}(a^2+b^2) } \Bigg)
\end{align*}
where we have used (\ref{h.def}), (\ref{h.1}) in the first inequality. The second inequality uses the fact that $\|E_i(t)\|=1$. Consider $r^* = \min ( \frac{\pi}{2\sqrt{\Delta}} , \frac{1}{\sqrt{\Delta}} \tan^{-1} (\frac{\sqrt{\Delta}(a^2+b^2)}{a+b})) $. One can easily verify that $r<r^*$, the RHS is positive in the above inequality.

\bibliographystyle{IEEEtran}
\bibliography{biblo}

 \end{document}